\RequirePackage{fix-cm}
\documentclass{svjour3}                     
\smartqed  
\usepackage{graphicx}
\usepackage{amsmath}
\usepackage{amssymb}
\usepackage{caption}
\usepackage{subcaption}
\usepackage{xcolor}
\usepackage{tablefootnote}

\begin{document}

\title{
On a workflow for efficient computation of the permeability of tight sandstones

\thanks{This work was supported by the Ministry of Science and Higher Education of the Russian Federation under agreement No. 075-10-2020-119 within the framework of the development program for a world-class Research Center.}
}



\authorrunning{Short form of author list} 

\institute{V. Pimanov \at
              Skolkovo Institute of Science and Technology,
                Bolshoy Boulevard 30, bld. 1, Moscow, Russia, 121205 \\
              \email{Vladislav.Pimanov@skoltech.ru}           
    \and
            V. Lukoshkin \at
        Skolkovo Institute of Science and Technology 
    \and
            P. Toktaliev \at
        Fraunhofer Institute for Industrial Mathematics
    \and
           O. Iliev \at
        Fraunhofer Institute for Industrial Mathematics
    \and
           E. Muravleva \at
        Skolkovo Institute of Science and Technology
    \and
           D. Orlov \at
        Skolkovo Institute of Science and Technology
    \and
           V. Krutko \at
        Gazpromneft Science and Technology Center
    \and
           A. Avdonin \at
        Gazpromneft Science and Technology Center
    \and
           K. Steiner \at
        Fraunhofer Institute for Industrial Mathematics
    \and
           D. Koroteev \at
        Skolkovo Institute of Science and Technology
}

\date{Received: date / Accepted: date}

\authorrunning{Pimanov et al.}
\maketitle

\newpage
\textcolor{blue}{\Large 
}

\tableofcontents

\newpage
\begin{abstract}
The paper presents a workflow for fast pore-scale simulation of single-phase flow in tight reservoirs typically characterized by low, multiscale porosity. Multiscale porosity implies that the computational domain contains porous voxels (unresolved porosity) in addition to pure fluid voxels. In this case, the Stokes-Brinkman equations govern the flow, with the Darcy term needed to account for the flow in the porous voxels. As the central part of our workflow, robust and efficient solvers for Stokes and Stokes-Brinkman equations are presented. The solvers are customized for low-porosity binary and multiclass images, respectively. Another essential component of the workflow is a preprocessing module for classifying images with respect to the connectivity of the multiscale pore space. 
{Particularly}, an approximation of the Stokes-Brinkman problem, namely, the Darcy problem, is investigated for the images that do not have pure fluid percolation paths. 
Thorough computational experiments demonstrate efficiency and robustness of the workflow for simulations on images from tight reservoirs.
Raw files describing the used CT images are provided as supplementary materials to enable other researchers to use them.

\end{abstract}

\section{Introduction.}
\label{Intro}

\subsection{Motivation and goals.}

Today the leading oil companies of the global energy perspective put forth great efforts to have commercial and stable production from tight and ultra-tight hydrocarbon resources in which the initial reserves are enormously more significant than the reserve of conventional resources. Despite all technological achievements, the development of unconventional resources might yet be a risky investment if it is not supported well with detailed and accurate reservoir studies \cite{wang2019review}. As a result and subject to the requirements of reservoir studies and relevant challenges, developing techniques, methods, and instruments for the comprehensive evaluation of unconventional hydrocarbon studies is currently an important area of research \cite{li2019determination}.  

Regardless of many attempts that have been made to modify the conventional methods of core analysis suitable to be applied for the evaluation of unconventional hydrocarbon resources \cite{lu2020analysis}, they are still time-consuming, expensive, and inaccurate \cite{sander2017laboratory}.  It is mainly due to the existence of sub-micron pores and throats that have strong effects on the storage and flow of tight and ultra-tight oil and gas resources \cite{chai2019new,javadpour2007nanoscale}.

To fit the risk attitude of the decision-makers, applying modern methods of core analysis like Digital Rock Physics, DRP, has mainly been focused on over the last decade \cite{kelly2016assessing,karoly2006psychological}. Besides the fact that DRP is an efficient cost control and risk management method, it allows the researchers to perform multiple numerical experiments on exactly the same sample and implement various analyses simultaneously \cite{jia2016some,verri2017development}.

Single-phase flows are attracting significant attention in DRP. A number of commercial \cite{thermofisher-web-page,geodict-web-page,openfoam-web-page} and academic \cite{Blunt2012,krotkiewski2011importance,koroteev2014direct,balashov2020dimp} flow solvers have been developed, to name just a few. These and many others are currently actively used for solving scientific, environmental, and industrial problems, as well as in benchmark studies \cite{andra2013digital,saxena2017references}. Despite the active development of algorithms and software in DRP, pore-scale simulations for tight reservoirs, usually characterized by low, multiscale porosity, remain a challenge. It should be noted that not only the tight reservoir rocks are characterized by multiscale porosity. There are also plenty of examples beyond natural rocks, e.g., active carbon particles in nonwoven filtering materials used for air filtration, catalyst in diesel particulate filters for automotive applications, bones, especially in the case of osteoporosis, etc. In such cases, despite the high imaging resolution, there is often unresolved porosity as a part of the image. Cutting-edge research is currently carried out to investigate the impact of the unresolved porosity on the considered problems. Due to the enormous complexity of the pore geometry, the variation of the material properties in the unresolved regions, the jumping coefficients, the large size of the computational domains, the existing numerical algorithms cannot always provide satisfactory results. The flow solvers developed for computing permeability for higher porosity in certain cases fail to converge for images of tight rocks. The demand for developing advanced customized algorithms in this area has risen essentially in the last years.

\subsection{Contributions.} 
The ultimate objective of this paper is to propose a workflow for solving single-phase pore-scale flow problems in tight sandstones characterized by low or very low porosity and connectivity, especially in the case of the appearance of unresolved porosity.
The workflow has essentially three components. \\ 
(i) The first stage, as usually, is the image processing. The multiclass models of real rock samples from tight reservoirs considered here are built according to the approach from \cite{avdonin2021application}, with the help of the double $\mu$CT scanning technique. The details of the model preparation process is a subject of follow up article, a sketch is briefly described in Section \ref{subsec:multiclass_preparation}. 
 \\
(ii) At the second stage, the 3D images (multiclass models) are classified with respect to the connectivity of the resolved pore space. \\
(iii) At the third stage, either the newly developed Stokes-Brinkman or Darcy solver is used to calculate the effective permeability for different classes of images. 


Extensive studies are carried out on the performance of the presented solvers. The simulation results illustrate the high accuracy and robustness of the Stokes-Brinkman solver. Its superior performance in computing permeability for tight rocks is demonstrated in comparison with state-of-the-art Stokes-Brinkman commercial solvers for DRP. 
The usage of the Darcy approximation of the Stokes-Brinkman equations to compute the permeability of images without pure fluid percolation paths is explored, and it shows a significant reduction of the CPU time, while providing good accuracy. 
Note that the idea of solving the Darcy problem instead of the Stokes-Brinkman one is not a new one; it was earlier explored, e.g., in \cite{krotkiewski2011importance}. Here this approach is a part of our workflow, which is equipped with a preprocessor module (classifier) to automatically decide if it is reasonable to solve the Darcy approximation.

Although the workflow is primarily intended for multiclass images, the robust Stokes solver for low porosity binary images is presented as well.


The remainder of the paper is organized as follows. Section \ref{sec:problem_statement} is devoted to the statement of the boundary value problem for single-phase flow in the case of resolved and unresolved porosity. Section \ref{sec:numerical_method} discusses the numerical algorithms and their implementation. Section \ref{sec:preprocessing} is dedicated to describing the preprocessing stage in the proposed workflow. Section \ref{sec:computational_experiments} is the core of the paper, it contains the results of the numerical experiments and their discussion. Data from the literature, when possible, and data from computations with the commercial software tool GeoDict \cite{geodict-web-page} are used for the validation. It starts with Subsection \ref{subsec:design_of_experiments} elaborating on the design of the computational experiments, describing the goals of the simulations. Subsection \ref{subsec:validation_binary} presents validation results for the Stokes solver on binary images. Subsection \ref{subsec:real_rocks} contains results of a thorough numerical investigation of the performance and accuracy of the developed Stokes-Brinkman and Darcy solvers on ternary and multiclass CT images from tight reservoir rocks. Subsection \ref{subsubsec:goals_and_samples} describes the collection of images used. In subsection \ref{subsubsec:validation_ternary}, simulation results for $300^3$ ternary samples are presented, and the applicability of the Darcy solver is discussed. Subsection \ref{subsubsec:multiclass} is dedicated to simulations with a multiclass image when each voxel has individual permeability. Subsection \ref{subsubsec:ternary_600_1350} is similar to subsection \ref{subsubsec:validation_ternary} but larger samples of sizes $600^3$ and $1350^3$ are considered there. Finally, in Subsection \ref{subsec:summary} some conclusions are drawn. Appendix containing a part of the results concludes the article.

\section{Problem statement.}\label{sec:problem_statement}

\subsection{Computational domain.}

We represent images by a unit cube domain $\overline{\Omega} = [0, 1]^3$ consisting of three non-overlapping parts:
\begin{equation}\label{domain_partition}
    {\Omega} = {\Omega}_f^h \cup {\Omega}_p^h \cup {\Omega}_s^h,
\end{equation}
where $\Omega_f^h$ is a pure fluid part (resolved porosity), $\Omega_s^h$ is a solid part, and $\Omega_p^h$ denotes a porous space (unresolved porosity). We suppose that these domains are voxelized 3D CT images, with images of tight sandstone being the main target. 
Namely, the domain ${\Omega}$ is composed from voxels:
$$\omega_{i,j,k} = [(i - 1)h; ih] \times [(j - 1)h; jh] \times [(k - 1)h; kh], \qquad i,j,k =  \overline{1, N}, \; h = 1/N,$$
where $N$ is the number of voxels in each dimension.
For each voxel $\omega_{i,j,k}$ we prescribe the porosity $\phi_{i,j,k} \in [0, 100] \subset \mathbb{N}$ which defines the partition \eqref{domain_partition} as follows:
\begin{equation}
\omega_{i,j,k} \subset
    \begin{cases}
     \overline{\Omega}_f^h \quad \mathrm{if} \quad \phi_{i,j,k} = 0, \\
     \overline{\Omega}_p^h \quad \mathrm{if} \quad \phi_{i,j,k} \in (0, 100), \\
     \overline{\Omega}_s^h \quad \mathrm{if} \quad \phi_{i,j,k} = 100.
    \end{cases}
\label{multiporosity}
\end{equation}
The solid part ${\Omega}^h_s$ is excluded from computations because no equations are prescribed there, and the computational domain becomes: 
\begin{equation*}
    \Omega^h = {\Omega}_f^h \cup {\Omega}_p^h.
\end{equation*}
Note that in \eqref{multiporosity} different porosity, coming from the grey image, can be considered in each porous voxel. We call such images {\it multiclass images}.
However, for the sake of simplicity, most simulations are performed for {\it ternary images} with some constant average porosity prescribed to all porous voxels.
Furthermore, {\it binary images} correspond to the case when all the porosity is resolved, and an image consists only of fluid and solid voxels (i.e., $\Omega^h = {\Omega}^h_f$).

\subsection{Governing equations.}

\noindent
{\bf Stokes problem for binary images.}
In the case when all the pores are resolved, the steady state slow flow of incompressible fluid is governed by the \emph{Stokes equations}:
\begin{equation}\label{stokes}
    \begin{split}
        -\Delta \mathbf{u} + \nabla p  &= \mathbf{0} \text{ in } \Omega^h, \\
        -\nabla \cdot \mathbf{u} &= 0 \text{ in } \Omega^h, \\
        + \; &\mathrm{b.c. \; on} \; \partial \Omega^h,
    \end{split}
\end{equation}
where $p$ is the fluid pressure and $\mathbf{u}$ is the fluid velocity. The boundary conditions on $\partial \Omega^h$ will be discussed in Section \ref{subsec:bc}.

\noindent
{\bf Stokes-Brinkman problem for ternary and multiclass images.}
There are two basic approaches for the numerical flow simulation when porous and free-flow regions co-exist at different scales.
The first approach is to use the coupled Stokes-Darcy equations (see, e.g., \cite{arbogast2006homogenization,layton2003darcystokes,discacciati2009darcystokes}) which consider of different models in different subdomains coupled via proper interface conditions. In DRP, the other approach, namely the Stokes-Brinkman model, sometimes also called the single domain approach, is valid in both the fluid and porous subregions (see, e.g., \cite{golfier2015investigation}). Note that often in the literature, such models are called just Brinkman models or Darcy-Brinkman models. However, \emph{Stokes-Brinkman} term is used here to emphasize that pure fluid voxels exist in the domain.
In $\Omega^h = {\Omega}^h_p \cup  {\Omega}^h_f$ we consider the \emph{Stokes-Brinkman equations}:
\begin{equation}\label{stokes_brinkman}
    \begin{split}
        -\Delta \mathbf{u} + \mathbf{K}^{-1} \mathbf{u} + \nabla p  &= \mathbf{0} \; \mathrm{in} \; \Omega^h, \quad \mathbf{K}^{-1} \geq 0,\\
        -\nabla \cdot \mathbf{u} &= 0 \; \mathrm{in} \; \Omega^h, \\
        + \; &\mathrm{b.c. \; on} \; \partial \Omega^h,
    \end{split}
\end{equation}
where $\mathbf{K}^{-1}$ is the dimensionless inverse permeability tensor, which is scaled by $L^2$ where $L$ [m] is a given characteristic (physical) length of a sample.
Note that, in the fluid part of the domain $\Omega^h_f$, 
$\mathbf{K}^{-1}$ becomes zero, and formally, one could recognize Stokes equations \eqref{stokes} in the resolved pores. 

\noindent
{\bf Brinkman perturbation of the Stokes-Brinkman problem. }
In tight reservoirs, one may face the necessity to compute the permeability for images where no pure fluid percolation path exists. In such cases, the flow consequently passes fluid and porous subregions; hence the effective permeability of the image is governed by the resistance of the porous subregions, while the resistance of the fluid subdomain could be neglected. For such samples, the perturbation $\mathbf{\tilde{K}}$ of the permeability tensor $\mathbf{K}$ can be considered by introducing a fictitious permeability value ${K}_{Stokes}$ in the pure fluid voxels from $\Omega_f^h$.
Then, $\mathbf{\tilde{K}}^{-1}$ does not vanish in any of the voxels in $\Omega^h$, and the \emph{Stokes-Brinkman equations} \eqref{stokes_brinkman} are thus perturbed to the \emph{Brinkman equations}:
\begin{equation}\label{brinkman}
    \begin{split}
        -\Delta \mathbf{u} + \mathbf{\tilde{K}}^{-1} \mathbf{u} + \nabla p  &= \mathbf{0} \; \mathrm{in} \; \Omega^h, \quad \mathbf{\tilde{K}}^{-1} > 0,\\
        -\nabla \cdot \mathbf{u} &= 0 \; \mathrm{in} \; \Omega^h, \\
        + \; &\mathrm{b.c. \; on} \; \partial \Omega^h.
    \end{split}
\end{equation}
The Brinkman equations (see, e.g.,  \cite{brinkman1949calculation}) play an important role for classes of porous media flows, e.g., describing more adequately flows in large porosity domains and allowing to account for no-slip boundary conditions. Here, however, they are not used for a more accurate description of the flow, but they are considered as a perturbation to the Stokes-Brinkman equations \eqref{stokes_brinkman}.

\noindent
{\bf Darcy approximation of the Brinkman problem. }
The Brinkman problem \eqref{brinkman} can be further approximated by neglecting the viscous term, which results in the \emph{Darcy equations}:

\begin{equation}\label{darcy}
    \begin{split}
        \mathbf{\tilde{K}}^{-1} \mathbf{u} + \nabla p  &= \mathbf{0} \text{ in } \Omega^h, \quad \mathbf{\tilde{K}}^{-1} > 0, \\
        -\nabla \cdot \mathbf{u} &= 0 \text{ in } \Omega^h, \\
        + \; &\mathrm{b.c. \; on} \; \partial \Omega^h.
    \end{split}
\end{equation}
It will be shown in this paper that for a certain class of images, the \emph{Darcy equations} \eqref{darcy} can be solved instead of the \emph{Stokes-Brinkman equations} \eqref{stokes_brinkman} in order to dramatically reduce the computational time while preserving the accuracy of the effective permeability computations.

\subsection{Boundary conditions.}\label{subsec:bc}

Assume, without loss of generality, that $\mathrm{z}$ direction is fixed as {\it flow direction}, and $\mathrm{x},\mathrm{y}$ directions as {\it tangential directions}. Then, the whole boundary $\partial \Omega^h$ is subdivided as follows:
\begin{equation*}
    \partial \Omega^h = \Gamma_{\mathrm{in}} \cup \Gamma_{\mathrm{out}} \cup \Gamma_\mathrm{t} \cup \Gamma_0,
\end{equation*}
  where  $\Gamma_{\mathrm{in}} = \{\mathbf{x} \in \partial \Omega^h \; | \; \mathrm{z} = 0\}$ is the inlet boundary, $\Gamma_{\mathrm{out}} = \{\mathbf{x} \in \partial \Omega^h \; | \; \mathrm{z} = 1\}$ is the outlet boundary, $\Gamma_{\mathrm{t}} = \{\mathbf{x} \in \partial \Omega^h \; | \; \mathrm{x} = 0,1 \text{ or } \mathrm{y} = 0,1\}$ denotes boundaries in the tangential directions, and $\Gamma_0 =  {\overline{\Omega}}^h_s \cap ( {\overline{\Omega}}^h_f \cup  {\overline{\Omega}}^h_p)$ denotes internal boundaries.
  
  The rigorous computation of the permeability tensor according to the homogenization theory is based on using periodic boundary conditions, see, e.g., \cite{hornung1997homogenization}. In the engineering practice, another formulation is more often used; namely, the velocity is prescribed at the inlet $\Gamma_{\mathrm{in}}$, and the pressure is prescribed at the outlet $\Gamma_{\mathrm{out}}$. It was shown in \cite{griebel2010homogenization} that the results obtained with this velocity/pressure formulation are equivalent to the results obtained with periodic boundary conditions. 
Furthermore, recently in the engineering literature, one more formulation was introduced, namely, the so-called \emph{pressure drop} boundary condition \cite{guibert2015computational}.
In this case, the following set of boundary conditions for equations \eqref{stokes}-\eqref{darcy} is prescribed:
\begin{equation}\label{pipo_bc}
    \begin{cases}
        p &= p_{\mathrm{in\phantom{p}}}, \quad \frac{\partial \mathbf{u}}{\partial \mathbf{n}} = \mathbf{0} \; \mathrm{on} \; \Gamma_{\mathrm{in}}, \\
        p &= p_{\mathrm{out}}, \quad \frac{\partial \mathbf{u}}{\partial \mathbf{n}} = \mathbf{0} \; \mathrm{on} \; \Gamma_{\mathrm{out}}, \\
        \mathbf{u} &= \mathbf{0} \; \mathrm{on} \; \Gamma_\mathrm{t}, \\
        \mathbf{u} &= \mathbf{0} \; \mathrm{on} \; \Gamma_0,
    \end{cases}
\end{equation}
where $(p_{\mathrm{out}} - p_{\mathrm{in}})$ [Pa] defines the pressure drop.
    
Two types of boundary conditions in flow direction are implemented in our workflow. First one is the pressure drop formulation \eqref{pipo_bc}, and the second one is the periodic formulation.

\subsection{Computing effective permeability.}

The effective permeability tensor, denoted $\mathbf{K^{eff}}$, is an intrinsic characteristic of the porous geometry, and it is a functional of the solution $\mathbf{u}$. Computation of the permeability tensor according to the homogenization theory can be found in \cite{hornung1997homogenization}. Here, the approach most often used in the engineering literature is considered. Orthotropic (diagonal) permeability tensor is assumed: $\mathbf{K^{eff}} = \mathrm{diag} (k_{\mathrm{xx}}^{eff}, k_{\mathrm{yy}}^{eff}, k_{\mathrm{zz}}^{eff})$. For a selected flow direction $\mathrm{z}$, and for a given pressure drop, the respective component of the permeability tensor $ k_{\mathrm{zz}}^{eff} [m ^ 2] $  is determined according to the Darcy's law:
\begin{equation}\label{darcy_law}
    k_{\mathrm{zz}}^{eff} = \hat{k}_{\mathrm{zz}}^{eff} L^2  = - \dfrac{\langle u_\mathrm{z} \rangle}{p_{\mathrm{out}} - p_{\mathrm{in}}} L^2,
\end{equation}
where $\hat{k}_{\mathrm{zz}}^{eff}$ is the corresponding component of the dimensionless permeability tensor, and
the Darcy's velocity $\langle u_\mathrm{z} \rangle$ is calculated by averaging the respective velocity component over the entire volume of the porous sample \cite{whitaker1986flow}:
\begin{equation}
    \langle u_\mathrm{z} \rangle = {|{\Omega}|}^{-1}\int_{\Omega^h}u_\mathrm{z} dV.
\end{equation}
It should be noted that in practice, to calculate effective permeability, the following formula is typically used:
\begin{equation}\label{darcy_law_flux}
    k_{\mathrm{zz}}^{eff} = -\dfrac{1}{p_{\mathrm{out}} - p_{\mathrm{in}}}\dfrac{Q}{A} L^2,
\end{equation}
where $ Q [m ^ 3 / s] $ is the volumetric flow rate, $ A [m ^ 2] $ is the cross-sectional area of the sample, and the Darcy velocity approximates the flow:
\begin{equation}
    \dfrac{Q}{A} \approx \langle u_\mathrm{z} \rangle.
\end{equation}


\section{Numerical method.}
\label{sec:numerical_method}


\subsection{Discretization.}
The discretization of the Stokes \eqref{stokes}, (Stokes-)Brinkman \eqref{stokes_brinkman}-\eqref{brinkman}, or Darcy \eqref{darcy} equations results in a block system of linear equations of the following form:
\begin{equation}\label{coupled_matrix}
    \begin{bmatrix}
    {A} & {B}^T \\ 
    {B} & 0
    \end{bmatrix}
    \begin{bmatrix}
    {u}_h \\ 
    p_h
    \end{bmatrix}
    =
    \begin{bmatrix}
    {f}_h \\ 
    0
    \end{bmatrix},
\end{equation}
where:
\begin{equation*}
A =
    \begin{cases}
     -{\Lambda}, &\quad \text{for equations \eqref{stokes}},\\
     -{\Lambda} + {K}^{-1}, &\quad \text{for equations \eqref{stokes_brinkman}-\eqref{brinkman}},\\
     \phantom{-} {K}^{-1}, &\quad \text{for equations \eqref{darcy}}.
    \end{cases}
\end{equation*}
Matrices ${\Lambda}$, ${B}$, and ${B}^T$ are discretisations of the velocity vector Laplacian, the pressure gradient, and the velocity divergence operators, respectively.
Matrix $K^{-1}$ denotes the inverse permeability matrix that may have zeros on its diagonal in the Stokes-Brinkman case. The nonzeros in the right-hand side $f_h$ appear due to imposing the Dirichlet boundary conditions on the pressure $p$ (or as a volumetric force in case of periodic boundary conditions).

\noindent
{\bf Iterative method. } 
In the Pressure Schur Complement approach (according to the terminology from \cite{turek1999efficient}), the velocity is eliminated from \eqref{coupled_matrix} so that an equivalent problem is obtained: 
\begin{equation}\label{pressure_equation}
    S p_h = g_h, \quad \mathrm{where} \quad S = BA^{-1}B^T, \quad g_h = BA^{-1}f_h.
\end{equation}
Once the pressure $p_h$ is computed, the velocity field is reconstructed:
\begin{equation*}
    u_h = A^{-1} (f_h - B^T p_h).
\end{equation*}

\noindent
The preconditioned conjugate gradient (PCG) method \cite{verfurth1984combined,elman1994inexact} is used to solve equation \eqref{pressure_equation}:
\begin{equation}\label{pcg_outer}
     p_h^{k+1}=p_h^k - \alpha_k \hat{S}^{-1} \left(S p_h^k - g_h \right), 
 \end{equation}
 where $\alpha_k$ is a parameter of the Krylov subspace method (see e.g. \cite{saad2003iterative}, Section 9.2).
As a preconditioner $\hat{S}$, we consider the following approximation of the Schur complement matrix $S$:
\begin{equation}\label{simple}
    \hat{S} = B D^{-1} B^T, \qquad D = \mathrm{diag}(A).
\end{equation}
Two-stage inner-outer iterative process is used to solve \eqref{pressure_equation} to avoid computations with full matrices $A^{-1}$ and $S$.
At each iterative step of the outer iterations given by \eqref{pcg_outer}, two inner iterative processes are needed. Firstly, the application of $S$ requires an inversion of the velocity matrix $A$, and, secondly, the preconditioner $\hat{S}$ \eqref{simple} needs to be inverted. For both cases, the PCG method with Algebraic Multigrid preconditioner is applied. It is also worth clarifying that for the degenerate Darcy case \eqref{darcy}, the preconditioner $\hat{S}$ coincides with $S$, and the outer iterations, as well as the inner iterations for computing the diagonal inverse $A^{-1}$, trivially converge in one step.

It is worth emphasizing, that the classical SIMPLE (Semi-Implicit algorithm for Pressure Linked Equations) algorithm \cite{patankar1983calculation,benzi2005numerical,keyes2013multiphysics} relies on 
the same approximation $\hat{S}$ \eqref{simple} of the Schur complement $S$. The connection between SIMPLE and $\hat{S}$ is explained, e.g. in \cite{turek1999efficient} (Section 2.2). 
Namely, the classical SIMPLE can be considered as Richardson method for solving equation \eqref{pressure_equation}  preconditioned with $\hat{S}$. 
In conjunction with pore-scale single-phase flow simulation SIMPLE is used, e.g., in \cite{Blunt2012,linden2018specialized,PoreChem2017}. 
SIMPLE was also used as a smoother in Geometric Multigrid when computing effective permeability of academic and real porous geometries, see \cite{iliev2001flexible}. Unlike those papers, here $\hat{S}$ is used as a preconditioner for the Conjugate Gradient method, which is known to converge much faster than Richardson method.

\subsection{Implementations details.}
The above algorithm is implemented in the framework SCoPeS, which includes modules SCoPeS-S, SCoPeS-SB, and SCoPeS-D for solving Stokes, Stokes-Brinkman, and Darcy problems, respectively.
The solver is built on top of the PETSc open-source library \cite{petsc-web-page,petsc-user-ref}, which provides a linear algebra backend for the scalable (MPI parallel) solution of partial differential equations, including data structures for sparse matrix computations, a context for Krylov subspace methods (PETSc.KSP module), and a variety of preconditioners (PETSc.PC module). 
Particularly, the PCG method that is used for inverting $S$, $\hat{S}$, and $A$ is implemented within PETSc.KSP.KSPCG routine.

\noindent
{\bf Staggered grid.} 
For discretization, we exploit the classical finite-difference scheme on staggered grids \cite{harlow1965numerical}.
It is worth mentioning that in DRP, the inverse permeability tensor $\mathbf{K}^{-1}$ ($\mathbf{\tilde{K}}^{-1}$) is a piece-wise constant function of the porosity map $\phi_{i,j,k}$, but the degrees of freedom for the velocity components on staggered grids are located on the voxel faces. Therefore, either porosity or permeability should be averaged over the neighbor voxels (arithmetical mean), and the matrix $K^{-1}$ should be constructed respectively.

\noindent
{\bf Solving with $A$ and $\hat{S}$.} 
When we solve \eqref{pressure_equation}, most of the time is spent on inverting matrices $A$ and $\hat{S}$ at each step of the outer iterations \eqref{pcg_outer}. The overall efficiency of the method is determined by the efficiency of the constructed preconditioners for $A$ and $\hat{S}$.
Multigrid methods are proven to work well for elliptic equations. However, the classical Geometric Multigrid method is not easy to apply in the case of complex pore-scale geometry. 
Another difficulty is related to strongly heterogeneous coefficients relevant for the matrix $\hat{S}$.
The Algebraic Multigrid approach generalizes the principles of Geometric Multigrid method for complex geometries and discontinuous coefficients \cite{ruge1987algebraic}.
We use BoomerAMG \cite{yang2002boomeramg}, which has demonstrated excellent performance in solving DRP problems.

\noindent
{\bf Stopping criteria.}
Let $rtol_S$, $rtol_A$, and $rtol_{\hat{S}}$ be the tolerances for inverting $S$, $A$, and $\hat{S}$, respectively.
The stopping criteria for the PCG iterations is determined by the relative residual in unpreconditioned norm
\newline
(PETSC.KSP\_NORM\_UNPRECONDITIONED). For example, the outer PCG iterative process \eqref{pcg_outer} stops as soon as:
\begin{equation}\label{unpreconditioned_norm}
    \dfrac{\|Sp_h^* - g_h\|_2}{\|g_h\|_2} \leq rtol_{S}.
\end{equation}
Given $rtol_S$ as an input of the algorithm, we determine $rtol_A$, and $rtol_{\hat{S}}$ as follows:
$$rtol_A = 10^{-2} \cdot rtol_{S}, \quad rtol_{\hat{S}} = rtol_{S}.$$



\section{Preprocessing stage.}
\label{sec:preprocessing}

\subsection{Multiclass model preparation.}
\label{subsec:multiclass_preparation}
Digital Rock Physics emerged recently as the alternative to conventional laboratory experiments on core samples. For  tight reservoir rocks, the usage of DRP is not as straight forward as for higher porosity sandstones. This is due to the inherent trade-off between the spatial resolution of data and the representativeness of the size of the model. For this regime, we developed a new approach to consider in a single 3D digital core model porosity from different scales (micro and sub-micro). 

To build a multiclass model of a core sample, we performed double X-ray $\mu$CT scanning. Before scanning, the sample was seated in the coreholder. The first scan was done for the sample fully saturated with air. Then after air evacuation, we saturated the sample in coreholder with Xenon (Xe) gas. Xe is a non-reactive and highly mobile gas that also has a lower radiolucency than grains. Comparison of two registered and calibrated 3D $\mu$CT data allows to map the distribution of Xe molecules and their amount in the volume of the digital core model. Assuming linear dependence between Xe intensity on the $\mu$CT data difference ($\mu$CT in air - $\mu$CT in Xe) and porosity of the 3D model voxels, it was possible to create a multiclass model. A specific porosity value characterized each model class with a discreteness level equal to 1\%. Multiclass model creation includes a few steps:
\begin{itemize}
    \item Data acquisition for the same sample position, the same resolution and the same X-ray parameters (electrical voltage and current, distance between X-ray source and sample centre, detector exposition and step of sample rotation).
    \item Image preprocessing (denoising, artefact removing).
    \item Registration of 3D images (matching images in space) .
    \item Intensity equilibration (to be sure that minerals without Xe have the same intensity on both 3D images).
    \item 3D images subtraction (to map only Xe distribution in model volume).
    \item Parametrization of dependence between Xe intensity and submicron porosity assuming that on subtracted. images, minimum (zero) intensity corresponds to 0\% porosity (solid) and maximum intensity corresponds to 100\% porosity (big pores).
    \item Applying determined “intensity-porosity” correlation to all 3D datasets of subtracted images.
\end{itemize}

It is also possible to construct a multiclass model by combining information about pore space structure from various sources of data: µCT, scanning microscopy (SEM), and focused ion beam with scanning microscopy (FIB-SEM). One can train a deep neural network to obtain a porosity map directly from SEM  or FIB-SEM images. Another opportunity is to create a multiclass model by coarsening a high-resolution binary model to reduce the size of a computational domain but remain the same physical scale of the model.

\subsection{Sample classification.}
A preprocessor is implemented as a separate module in C++ language; it relies on the Disjoint Set Union data structure \cite{galler1964improved}. First of all, as it is usual in DRP, it identifies isolated fluid (pure fluid or/and porous) subdomains and removes them. Next, the image is classified. 
The existence of percolation patch(es) in the computational domain $\Omega^h = \Omega^h_p \cup \Omega^h_f$ is checked. If percolation exists, the Stokes-Brinkman problem \eqref{stokes_brinkman} is well-posed in the computational domain $\Omega_h$, and the image is classified as an image of Category A. Furthermore, the connectivity of the pure fluid region $\Omega^h_f$ is checked. If percolation exists, the image is moved from Category A to Category B, images with Stokes connectivity. Note that after such a preprocessing, Category A contains images that do not have Stokes connectivity but have Stokes-Brinkman connectivity via sequences of pure fluid and porous voxels. 





\section{Computational experiments.}
\label{sec:computational_experiments}

\subsection{Design of computational experiments.}
\label{subsec:design_of_experiments}

Computational experiments have to be performed for the following purposes: 

(i) The accuracy of the developed solvers has to be validated in comparison with data from the literature or in comparison with simulations performed with validated solvers. 3D CT images and generated geometries can be used here. 

(ii) The impact of the stopping criteria, the usage of different boundary conditions, and robustness with respect to the problem size have to be studied.

(iii) The main parts of the computer simulations have to be carried out to evaluate the performance of the developed customized solvers for pore-scale simulation on 3D CT images of tight sandstones in the case of unresolved porosity, what is the main goal of this paper.


\subsection{Stokes solver. Validation and performance studies on binary images.}
\label{subsec:validation_binary}

\begin{figure}
     \centering
     \begin{subfigure}[b]{0.5\textwidth}
         \includegraphics[width=\textwidth]{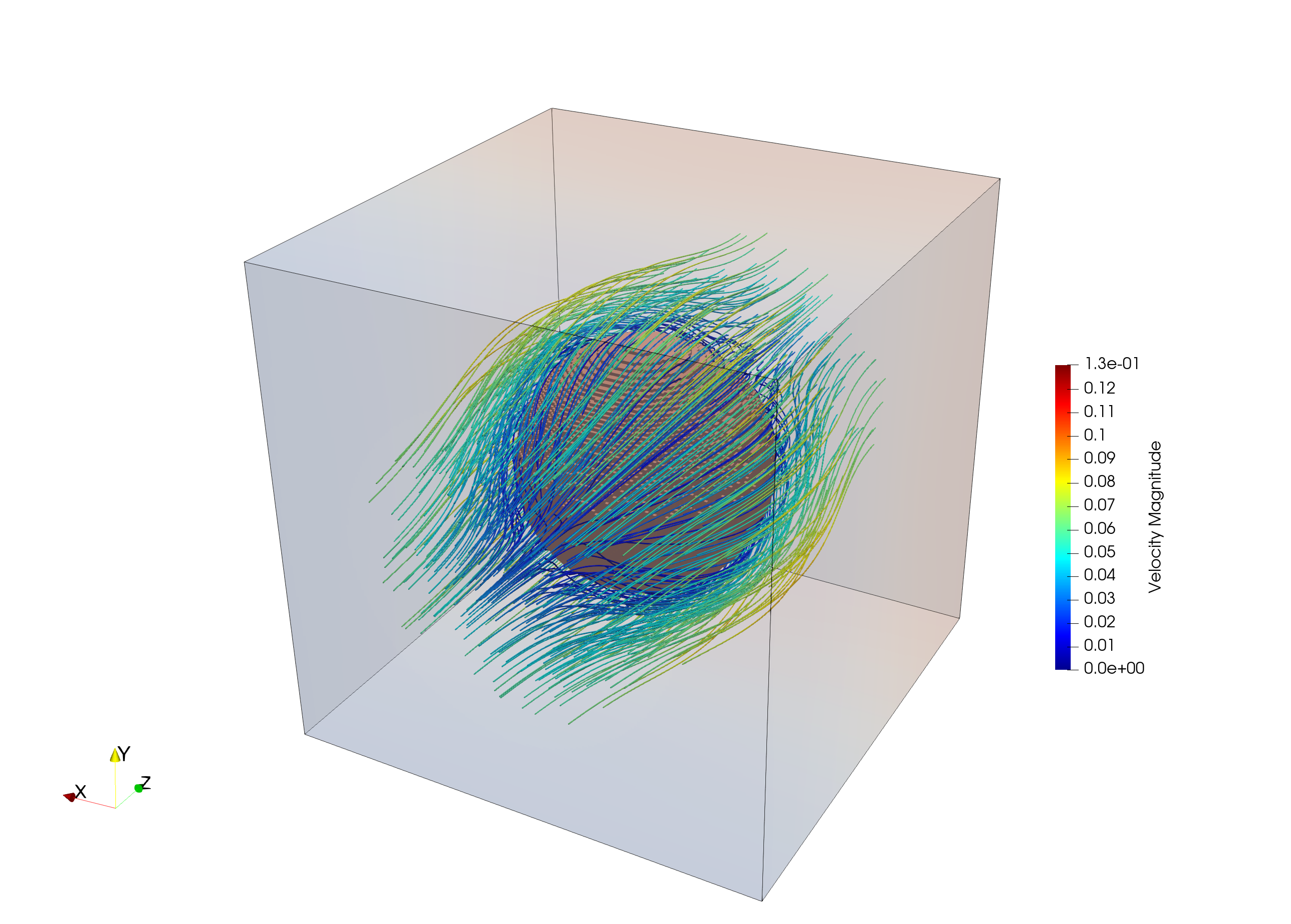}
         \caption{}
         \label{fig:sphere}
     \end{subfigure}
     \begin{subfigure}[b]{0.35\textwidth}
         \includegraphics[width=\textwidth]{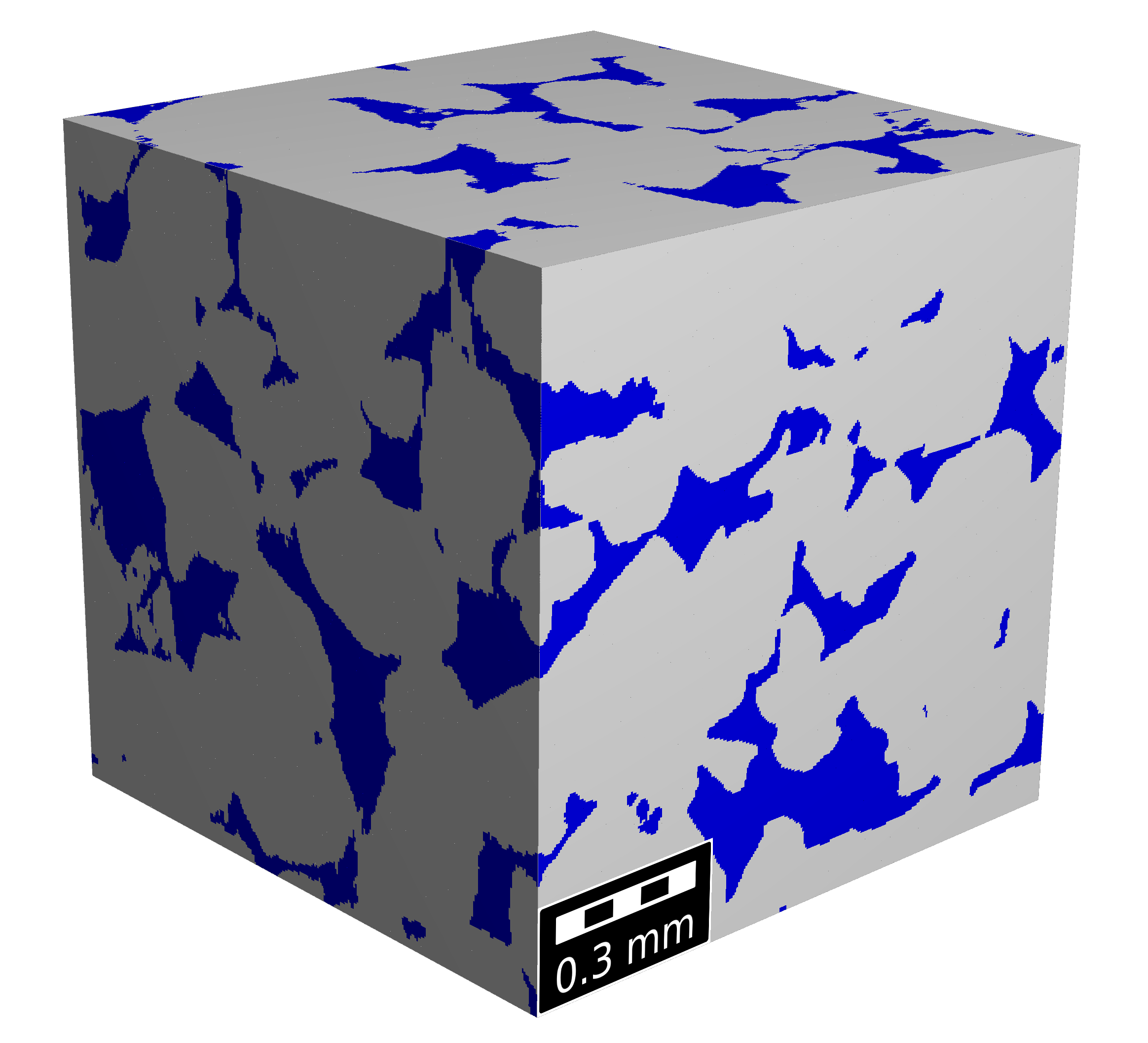}
         \caption{}
         \label{fig:binary300}
     \end{subfigure}
        \caption{(a) Velocity magnitude streamlines for periodic arrangements of spheres, diameter $D=0.5$, resolution $N = 160$. (b) Binary image of size $300^3$ with moderate-low porosity $\Phi = 0.21$, colors: grey - solid region, blue - free pores (GeoDict visualization).}
        \label{fig:sphere_binary300}
\end{figure}

{\bf Validation on periodic array of impermeable spheres.}
The developed Stokes solver, SCoPeS-S, is validated first on the classic example of flow around periodic arrangements of solid spheres. Imposing periodic boundary conditions in all three directions, one can consider only a single sphere. The computed permeability is compared with numerical and analytical results from  \cite{jung2005fluid} and \cite{sangani1982slow}. An exemplary computation is visualized on Figure \ref{fig:sphere}. The results are summarized in Table \ref{table:spheres}. The simulation results show that SCoPeS-S correctly computes the permeability. Also, convergence with respect to the grid size is clearly observed.





\begin{table}[h!]
\centering
 \begin{tabular}{||c c c c c c||} 
 \hline
 \hline
 D & N=40 & N=80 & N=160 & J\&T (ref \cite{jung2005fluid}) & S\&A (ref \cite{jung2005fluid}, \cite{sangani1982slow}) \\ [0.5ex] 
 \hline\hline
 0.1 & $9.74\cdot 10^{-1}$ & $9.01\cdot 10^{-1}$ & $9.02\cdot 10^{-1}$ &$9.15\cdot 10^{-1}$ & $9.11\cdot 10^{-1}$\\ 
 \hline
 0.2 & $3.77\cdot 10^{-1}$ & $3.78\cdot 10^{-1}$& $3.80\cdot 10^{-1}$ & $3.84\cdot 10^{-1}$ & $3.82\cdot 10^{-1}$\\
 \hline
 0.4 & $1.21\cdot 10^{-1}$ & $1.22\cdot 10^{-1}$ & $1.23\cdot 10^{-1}$ & $1.25\cdot 10^{-1}$ &$1.23\cdot 10^{-1}$\\
 \hline
 0.6 & $4.44\cdot 10^{-2}$ & $4.43\cdot 10^{-2}$ &  $4.43\cdot 10^{-2}$ & $4.58\cdot 10^{-2}$ & $4.45\cdot 10^{-2}$\\
 \hline
 0.8 & $1.29\cdot 10^{-2}$ & $1.31\cdot 10^{-2}$ & $1.31\cdot 10^{-2}$ & $1.38\cdot 10^{-2}$ &$1.32\cdot 10^{-2}$\\
 \hline
  1.0 & $2.48\cdot 10^{-3}$ & $2.51\cdot 10^{-3}$ & $2.51\cdot 10^{-3}$ & $2.67\cdot 10^{-3}$ & $2.52\cdot 10^{-3}$\\ [1ex] 
 \hline
\end{tabular}
\caption{Dimensionless permeability $\hat{k}_{\mathrm{zz}}^{eff}$ for periodic array of impermeable spheres. $D$ and $N$ denote the diameter of the spheres (with respect to the unit length of the cube) and the number of voxels in one direction, respectively. Relative tolerance $rtol_S = 10^{-3}$ was used in SCoPeS-S.}
\label{table:spheres}
\end{table}

\noindent
{\bf Validation on low porosity 3D CT images from the literature.} 
Further, validation of SCoPeS-S is done on five low porosity images listed in Table \ref{table:BinaryValidationLiterature}. The computed permeabilities are in good agreement with the data from the literature. The permeabilities computed with the software tools GeoDict, DiMP, and DHD are taken from \cite{orlov2021different}. 
It can be seen that the values computed with the newly developed solver, SCoPeS-S, are in the same range as the data from the literature.



\begin{table}[h!]
\centering
 \begin{tabular}{|c c c c c c|} 
 \hline
   & Porosity, $\Phi$ & GeoDict  & SCoPeS-S & DHD & DiMP\\
   \hline
 A & 0.058  & $7.1 \cdot 10^{-4}$ & $6.5 \cdot 10^{-4}$ & $6.7 \cdot 10^{-4}$ & $8.4 \cdot 10^{-4}$ \\
 B & 0.044  & $8.8 \cdot 10^{-4}$ & $8.1 \cdot 10^{-4}$ & $1.0 \cdot 10^{-3}$ & $8.7 \cdot 10^{-4}$ \\
 C & 0.053  & $3.9 \cdot 10^{-4}$ & $3.1 \cdot 10^{-4}$ & $3.2 \cdot 10^{-4}$ & $6.6 \cdot 10^{-4}$ \\
 D & 0.097  & $1.18 \cdot 10^{-2}$ & $9.83 \cdot 10^{-3}$ & $1.02 \cdot 10^{-2}$ & $9.38 \cdot 10^{-3}$ \\
 E & 0.064  & $8.3 \cdot 10^{-4}$ & $7.3 \cdot 10^{-4}$ & $6.8 \cdot 10^{-4}$ & $1.26 \cdot 10^{-3}$ \\
 \hline
\end{tabular}
\caption{Validation of SCoPeS-S on five low porosity binary images from \cite{orlov2021different}. Permeability $k_{\mathrm{zz}}^{eff}$ in $mkDa$, $L = 0.00072$ m, relative tolerance $rtol_S = 10^{-5}$.}
\label{table:BinaryValidationLiterature}
\end{table}




\noindent
{\bf Validation on 3D CT rock image with low-moderate porosity.} 
In order to further validate the Stokes solver, SCoPeS-S, cross validation with the Stokes solvers from GeoDict for binary $300^3$ image with porosity $\Phi=0.21$ was performed. 


\begin{table}[h!]
  \centering
  \begin{tabular}{||c c | c c c c ||}
    \hline
    \hline
     & SCoPeS & & \multicolumn{3}{c||}{GeoDict} \\
     & SCoPeS-S & & SimpleFFT & SimpleFFT  & LIR \\
     $rtol_S$ & Pin/Pout & Tol & Vin/Pout & Periodic & Periodic \\
    \hline
    \hline
    $ 10^{-2}$ & $6.74\cdot10^{6}$ (491) &$ 10^{-2}$ & $6.84\cdot10^{6}$ (3928) & $6.37\cdot10^{6}$ (1460) & $6.43\cdot10^{6}$ (325) \\
    \hline
    $10^{-3}$ & $6.45\cdot10^{6}$ (708) & $ 10^{-3}$ & $6.88\cdot10^{6}$ (6507) & $6.38\cdot10^{6}$ (2583) & $6.42\cdot10^{6}$ (771) \\
    \hline
    $10^{-4}$ & $6.46\cdot10^{6}$ (949) & $ 10^{-4}$ & $6.88\cdot10^{6}$ (7510) & $6.39\cdot10^{6}$ (3102) & $6.42\cdot10^{6}$ (1330) \\
    \hline
    \hline
  \end{tabular}
  \caption{SCoPeS-S and GeoDict (LIR, SimpleFFT) results for binary sample of size $300^3$ (Fig. \ref{fig:binary300}). Permeability $k_{\mathrm{zz}}^{eff}$ in $mkDa$ (CPU time in s). SCoPeS-S BC: pressure drop 1 Pa, $L=0.0011991$ m.}
  \label{table:BinaryImage}
\end{table}

\noindent
The results are summarized in Table \ref{table:BinaryImage}. It can be seen that SCoPeS-S provides results that are within the range of the computations done with GeoDict.  For each particular set of the computations, it was decided to gradually decrease the tolerance until the third digit in the computed permeability value is established. The efficiency of the SCoPeS Stokes solver is comparable to the LIR solver from GeoDict \cite{linden2015lir} but outperforms the solver SIMPLE-FFT \cite{wiegmann2007computation}.

\subsection{Stokes-Brinkman and Darcy solvers. Validation and performance studies on ternary and multiclass images from tight reservoir rocks.}
\label{subsec:real_rocks}

In this Section, the main simulation results are presented. The performance of the developed Stokes-Brinkman solver is systematically investigated for samples coming from tight reservoirs. The performance of SCoPeS-SB and GeoDict Stokes-Brinkman solvers is compared for the considered class of problems. Further on, as discussed above, just Darcy's problem is solved for images with no Stokes connectivity using SCoPeS-D solver.

\subsubsection{Preliminaries: Goals and samples.}
\label{subsubsec:goals_and_samples}

We start from three different multiclass samples of size $300^3$, named {\it S1,S2,S3}, which are sub-samples of a large multiclass  sample. Samples {\it S1,S2,S3} have $3.79, 3.25, 3.78$ mln. non-solid voxels, correspondingly. Recall that  {\it ternary images} are the images composed of solid voxels, fluid voxels, and identical porous voxels, the latter being responsible for all unresolved porosity.
Stokes-Brinkman equations have to be solved to compute the permeability for ternary images. For each of these three multiclass samples, the following procedure for increasing porosity and segmenting into ternary images is applied.
First, given a threshold $T \in (0, 100]$, we replace all porous voxels $\omega_{i,j,k} \in \Omega^h_p$ having porosity $\phi_{i,j,k} \geq T$ by pure fluid voxels. Next, the remaining porosity is arithmetically averaged with resulting value $\tilde{\phi}  \in (0, 100)$.
In this way, nine ternary samples are created, which are encoded by N\_T\_$\tilde{\phi}$.raw, where N=S1,S2,S3, T = 100,90,80, and $\tilde{\phi}$ is resulting averaged porosity. Note, $T=100$ corresponds to the case when the porosity is just averaged. 
Finally, constant permeability value is calculated for the averaged porosity $\tilde{\phi}$ according to \eqref{correlation} and assigned to all remaining porous voxels. 
 Table \ref{table:SamplesS123} summarizes the connectivity for each of the nine samples. Six of them are from Category A, and three of them are from Category B (marked in \textbf{bold}).
 
 The following correlation formula between porosity $\phi_{i,j,k} \in (0,100)$ and permeability $K_{i,j,k}$ [mkDa] was used: 
\begin{equation}\label{correlation}
    K_{i,j,k} = 7.251 \cdot 10^{-2} \exp(0.147076689 \cdot \phi_{i,j,k}).
\end{equation}
 This correlation formula \eqref{correlation} was derived using empirical petrophysical correlation from FIB-SEM data analysis \cite{avdonin2021application}.

 \begin{table}[h!]
	\centering
	\begin{tabular}{||c c c c c||} 
		\hline
		\hline
		Samples & $\Phi$ (Res, Unres) & Stokes  & Stokes-Brinkman & File name \\ [0.5ex] 
	  &   &   Connectivity &   Connectivity & \\ [0.5ex] 
		\hline
		\hline
		{S1\_100\_61} & 0.14 (0.044, 0.096) & No & Yes & S1\_100\_61.raw \\
		\hline
		{S1\_90\_56} & 0.14 (0.057, 0.083) & No & Yes & S1\_90\_56.raw \\
		\hline
		{S1\_80\_50} & 0.14 (0.070, 0.070) & No & Yes & S1\_80\_50.raw \\
		\hline
    		{S2\_100\_60} & 0.14 (0.050, 0.090) & No & Yes & S2\_100\_60.raw \\
		\hline
		{\bf S2\_90\_55} & { 0.14 (0.062, 0.078)} & {\bf Yes} & { Yes} & S2\_90\_55.raw\\
		\hline
		{\bf S2\_80\_49} & { 0.14 (0.075, 0.065)} & {\bf Yes} & { Yes} & S2\_80\_49.raw \\
		\hline
		{S3\_100\_58} & 0.12 (0.050, 0.070) & No & Yes & S3\_100\_58.raw \\
		\hline
		{S3\_90\_53} & 0.12 (0.057, 0.063) & No & Yes & S3\_90\_53.raw \\
		\hline
		{\bf S3\_80\_48} & { 0.12 (0.066, 0.054)} & {\bf Yes} & { Yes} & S3\_80\_48.raw \\
		\hline
	\end{tabular}
	\caption{Porosity and Connectivity of samples}
	\label{table:SamplesS123}
\end{table}
\noindent


The same segmentation procedure was also applied for two multiclass images of size $600^3$, named {\it T1,T2}, and having $30.5$, $30.7$ mln. non-solid voxels, respectively. Additionally, one large multiclass sample, {\it U1}, of size $1350^3$ with $350.8$ mln. non-solid voxels was studied. Porosities and connectivities of these samples are summarized in Table \ref{table:600connectivity}.  

This collection of images gives us a possibility to perform detailed testing of the performance of the Stokes-Brinkman solver for images with different porosity and different fraction of unresolved regions. The effective permeabilities computed with SCoPeS-SB solver is compared to the results obtained with GeoDict.

\begin{figure}
     \centering
     \begin{subfigure}[b]{0.4\textwidth}
         \includegraphics[width=\textwidth]{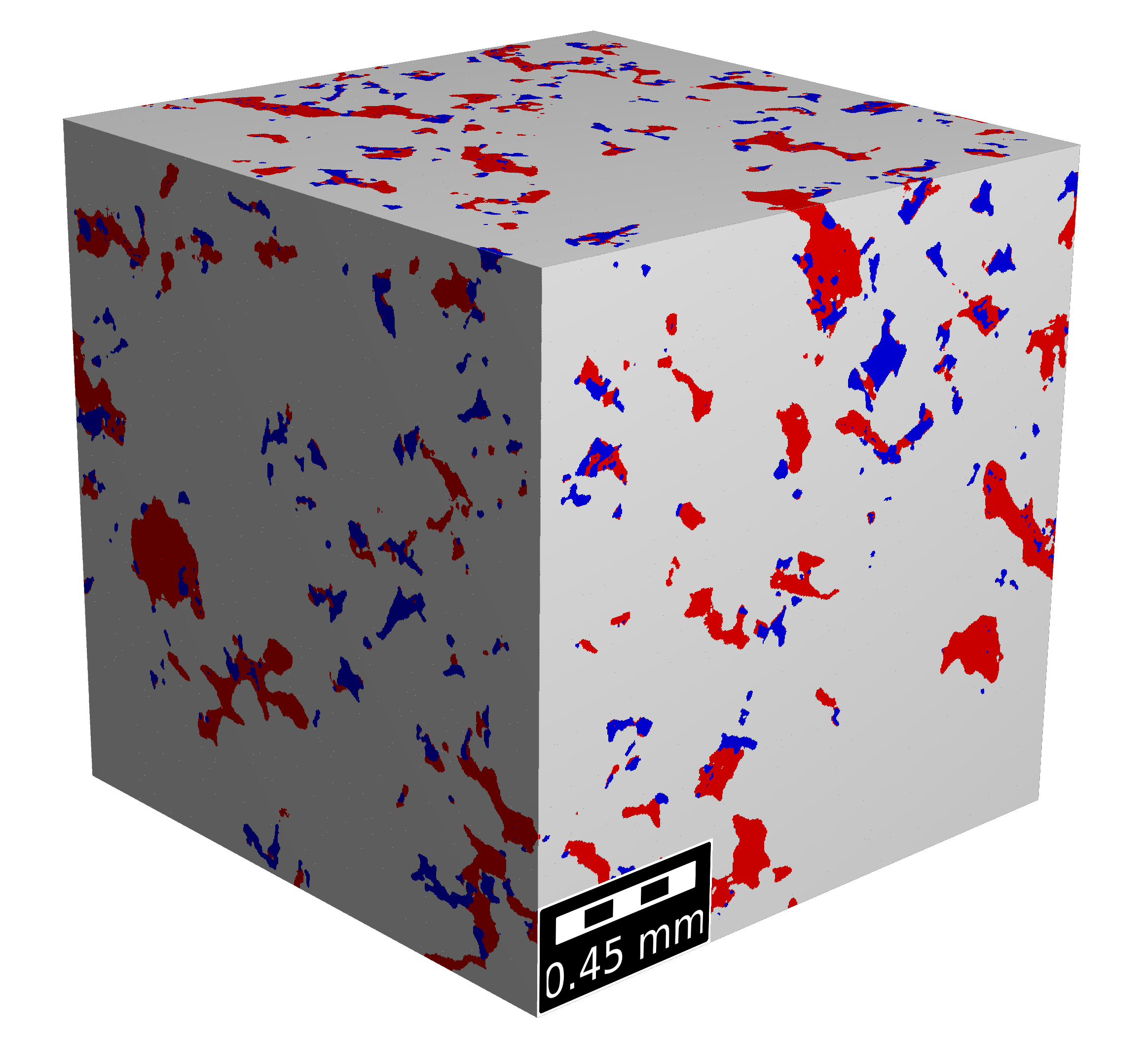}
         \caption{sample T1\_100\_59}
         \label{fig:T11}
     \end{subfigure}
     \begin{subfigure}[b]{0.4\textwidth}
         \includegraphics[width=\textwidth]{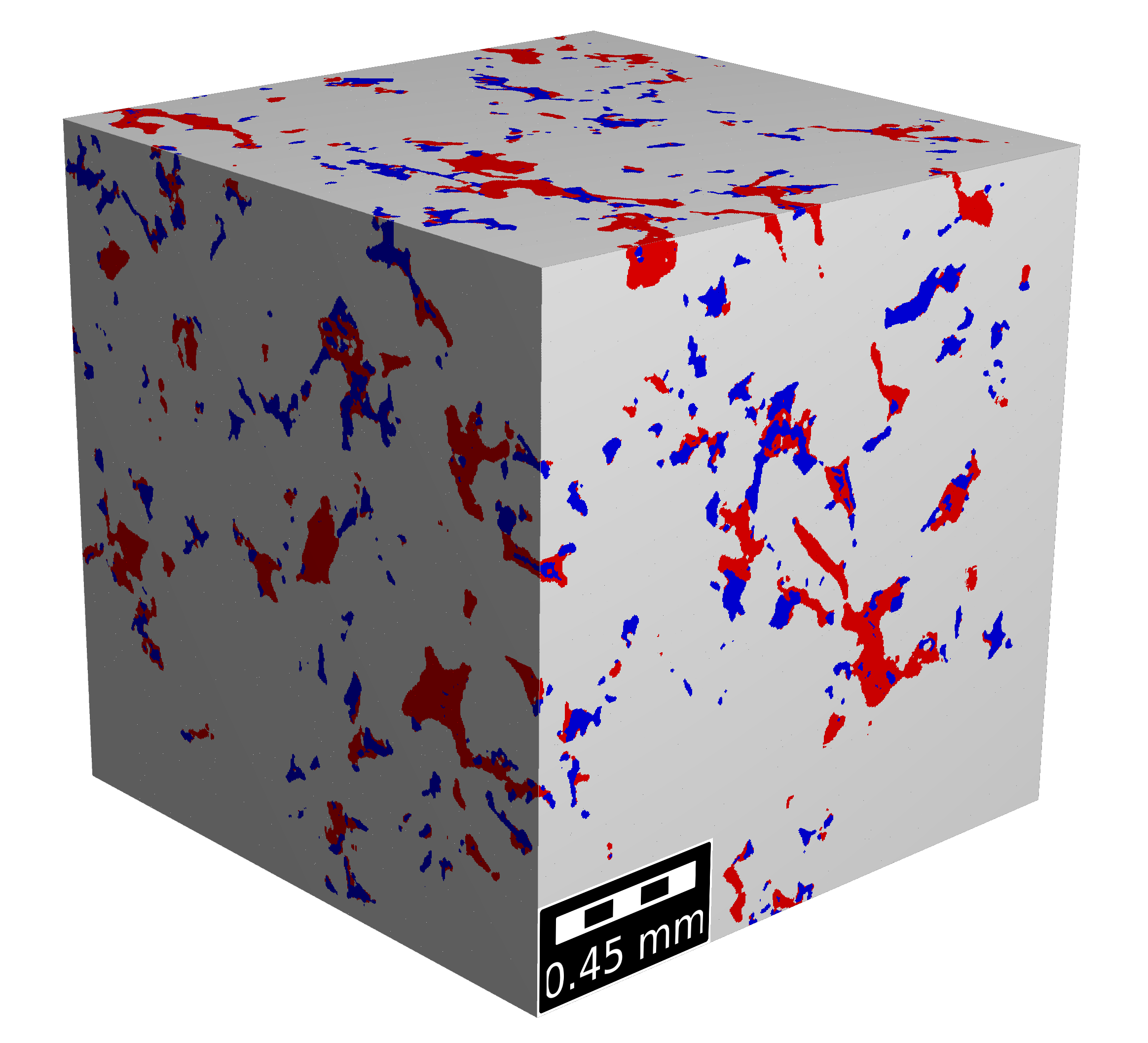}
         \caption{sample T2\_100\_58}
         \label{fig:T21}
     \end{subfigure}
        \caption{Ternary samples T1\_100\_59 and T2\_100\_58. Colors: grey - solid region, red - unresolved porosity region, blue - free pores. GeoDict visualization.}
        \label{fig:T1_T2_visualization}
\end{figure}

\begin{figure}[h]
\centering
\includegraphics[width=0.65\textwidth]{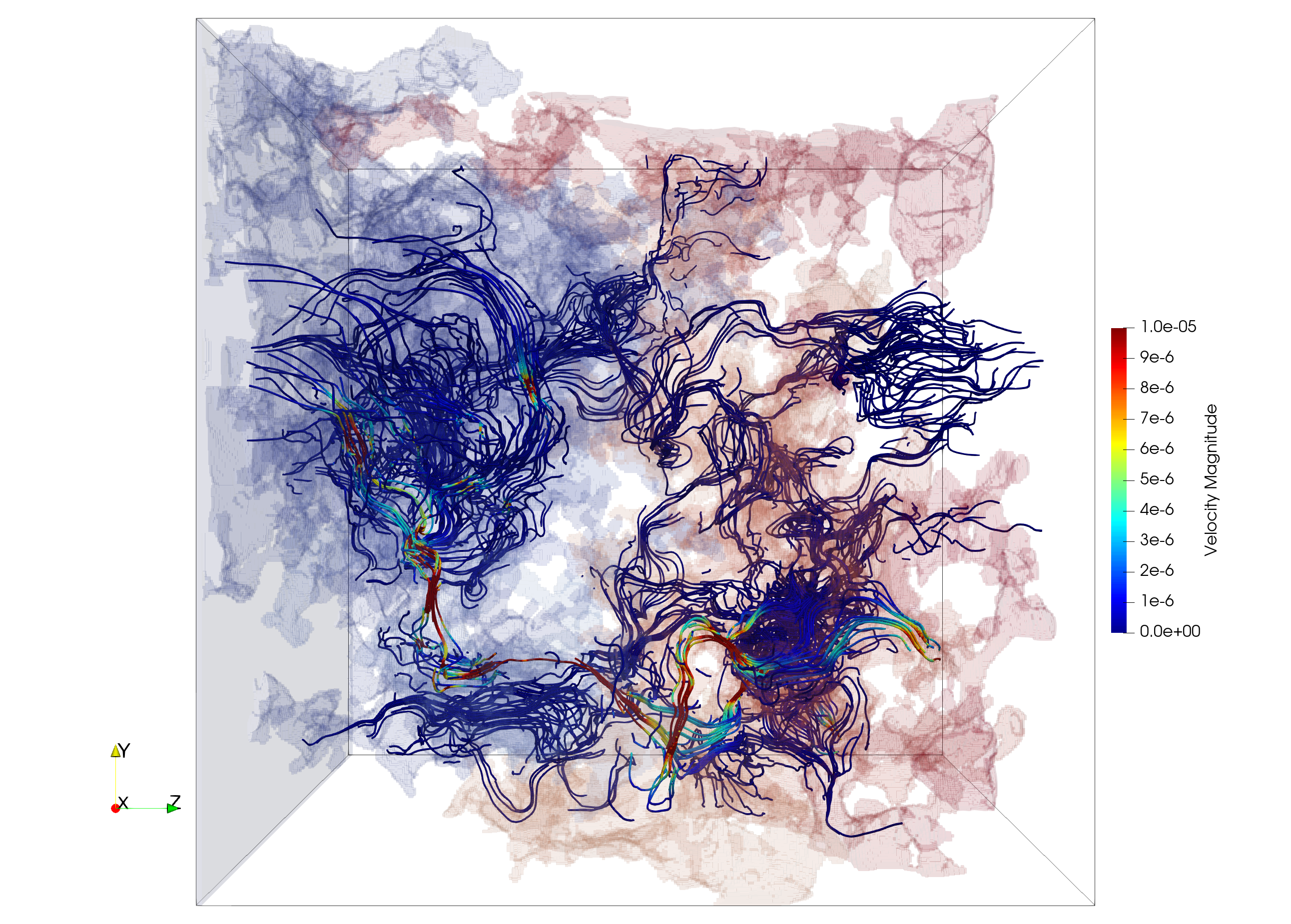}
\caption{Velocity magnitude streamlines. Paraview visualization for ternary sample S2\_80\_49, Stokes-Brinkman problem.}
\label{fig:s2_streamlines}
\end{figure}

Figure \ref{fig:T1_T2_visualization} presents explanatory pore-space visualization for two ternary samples: $T1\_100\_59$ and $T2\_100\_58$. Also, the result of Stokes-Brinkman simulation on ternary sample S2\_80\_49 is shown on Figure \ref{fig:s2_streamlines}.

\subsubsection{Simulations on ternary images of size $300^3$.}
\label{subsubsec:validation_ternary}

\noindent{\bf Solving Stokes-Brinkman equations for ternary images of Category A (no Stokes connectivity).}
Firstly, we consider three Category A images obtained from samples {\it S2} and {\it S3}. The simulation results of comparing the SCoPeS-SB solver and the SimpleFFT solver from GeoDict are presented in Table \ref{table:S2S3-A}. Also, the simulations results with the LIR solver from GeoDict are shown. However, LIR demonstrated poor performance for low-porosity images, therefore very few results were presented.
 It can be seen that SimpleFFT and LIR solvers converge very slowly but have the correct tendency, and the computed values are not far from what was expected. 
\begin{table}[h!]
	\centering
	\begin{tabular}{||c c c r | c c ||}
		\hline
		\multicolumn{6}{|| c ||}{Sample S2\_100\_60, Perm of porous voxels 493.0 $mkDa$} \\
		\hline
		\multicolumn{4}{|| c |}{GeoDict} & \multicolumn{2}{ c ||}{SCoPeS-SB} \\
		\hline
		Tol & Solver & $k_{\mathrm{zz}}^{eff}$, $mkDa$ & CPU, s & $rtol_S$ & $k_{\mathrm{zz}}^{eff}$, $mkDa$ (CPU, s) \\ [0.5ex]
		\hline
		\hline
		$ 10^{-1}$ & SimpleFFT: & $1.20 \cdot 10^2$ & 20511 & $ 10^{-6}$ & $6.80 \cdot 10^1$ (1564) \\
		\hline
		$ 1.2 \cdot 10^{-2}$ & SimpleFFT: & $7.95 \cdot 10^1$ & 172009 & $10^{-7}$ & $7.31 \cdot 10^1$ (2769) \\
		\hline
		 $ 1.9 \cdot 10^{-1}$ & LIR: & $1.04 \cdot 10^2$ & 214560 & $10^{-8}$ & $7.27 \cdot 10^1$ (3499) \\
		\hline
		 & & & & $ 10^{-9}$ & $7.27 \cdot 10^1$ (4090) \\
		\hline
		\hline
		\hline
		\multicolumn{6}{|| c ||}{Sample S3\_100\_58, Perm of porous voxels 367.4 $mkDa$} \\
		\hline
		\multicolumn{4}{|| c |}{GeoDict} & \multicolumn{2}{ c ||}{SCoPeS-SB} \\
		\hline
		Tol & Solver: & $k_{\mathrm{zz}}^{eff}$, $mkDa$ & CPU, s & $rtol_S$ & $k_{\mathrm{zz}}^{eff}$, $mkDa$ (CPU, s) \\ [0.5ex]
		\hline
		\hline
		$ 10^{-1}$ & SimpleFFT: & $4.07 \cdot 10^1$ & 5862 & $ 10^{-7}$ & $2.23 \cdot 10^1$ (1829) \\
		\hline
		$ 1.1 \cdot 10^{-2}$ & SimpleFFT: & $2.86 \cdot 10^1$ & 81435 & $ 10^{-8}$ & $2.54 \cdot 10^1$ (2199) \\
		\hline
		 & & & & $ 10^{-9}$ & $2.57 \cdot 10^1$ (2693) \\
		\hline
		 & & & & $ 10^{-10}$ & $2.57 \cdot 10^1$ (3265) \\
		\hline
		\hline
		\hline
		\multicolumn{6}{|| c ||}{Sample S3\_90\_53, Perm of porous voxels 176.0 $mkDa$} \\
		\hline
		\multicolumn{4}{|| c |}{GeoDict} & \multicolumn{2}{ c ||}{SCoPeS-SB} \\
		\hline
		Tol & Solver: & $k_{\mathrm{zz}}^{eff}$, $mkDa$ & CPU, s & $rtol_S$ & $k_{\mathrm{zz}}^{eff}$, $mkDa$ (CPU, s) \\ [0.5ex]
		\hline
		\hline
		$10^{-1}$ & SimpleFFT: & $6.52 \cdot 10^1$ & 10758 & $10^{-6}$ & $2.50 \cdot 10^1$ (1487) \\
		\hline
		$8\cdot10^{-2}$ & SimpleFFT: & $5.18 \cdot 10^1$ & 104421 & $10^{-7}$ & $5.12 \cdot 10^1$ (1875) \\
		\hline
		 & & & & $10^{-8}$ & $4.92 \cdot 10^1$ (2220) \\
		\hline
		 & & & & $10^{-9}$ & $4.92 \cdot 10^1$ (2733) \\
		\hline
		\hline
	\end{tabular}
	\caption{Ternary samples S2 and S3 of Category A. Permeability $k_{\mathrm{zz}}^{eff}$ in $mkDa$ computed with GeoDict (solvers SimpleFFT and LIR with periodic bc), and with SCoPeS-SB with pressure drop bc, $L=0.0009$ m, nproc=8.}
	\label{table:S2S3-A}
\end{table}

\noindent
Note, all three images produced from {\it S1} sample also belong to Category A, and the corresponding results are available in Table \ref{table:S1-A} in the Appendix.
In all cases, simulations with SCoPeS-SB show very robust convergence with respect to the selected tolerance. Unlike the GeoDict solvers, the computational time increases moderately when decreasing the tolerance value.

\noindent{\bf Solving Darcy approximation for ternary images of Category A (no Stokes connectivity).}
As it was mentioned above, in the case of no Stokes connectivity, it makes sense to first explore approximating Stokes-Brinkman equations with the Darcy equation, and after that solve it to compute the flow, and thus the permeability. The approximation is done by adding artificial permeability in the pure fluid voxels, and after that dropping the viscous terms. Substituting the velocity (in this case with a diagonal matrix) into the continuity equation, we obtain a scalar second order elliptic equation for the pressure. 

Simulation results with Darcy model for samples {\it S2} and {\it S3} are presented in Table \ref{table:Darcy-S2S3-A}.
The results are computed using Darcy solver SCoPeS-D. For comparison, the effective permeabilities computed with SCoPeS-SB solver are presented in the last line of the Table. One can see that the computations with the Darcy approximation are about 150 times faster compared to SCoPeS-SB, and at the same time the same accuracy in the computation of the permeability of the sample can be achieved.
Similar results for samples {\it S1} are available in Table \ref{table:Darcy-S1-A} in the Appendix.

Additionally, simulation results from sensitivity study on how the artificial permeability in the Darcy approximation influences the accuracy of the computations are presented in Table \ref{table:Darcy_S1} in the Appendix. The sample S1\_100\_61 is considered there.

\begin{table}[h!]
	\centering
	\begin{tabular}{||c c c c||}
		\hline
		\hline
		$K_{Stokes}$, $mkDa$ & S2\_100\_60 & S3\_100\_58 & S3\_90\_53 \\ [0.5ex]
		\hline
		\hline
		$10^5$ & $4.11 \cdot 10^1$ (42.8) & $2.34 \cdot 10^1$ (43.4) & $3.53 \cdot 10^1$ (42.6) \\
		\hline
		$ 10^7$ & $7.43 \cdot 10^1$ (44.1) & $2.57 \cdot 10^1$ (43.1) & $4.95 \cdot 10^1$ (44.2) \\
		\hline
		$10^9$ & $7.54 \cdot 10^1$ (44.8) & $2.58 \cdot 10^1$ (43.4) & $4.99 \cdot 10^1$ (43.5) \\
		\hline
		$ 10^{10}$ & $7.54 \cdot 10^1$ (44.5) & $2.58 \cdot 10^1$ (44.0) & $4.99 \cdot 10^1$ (43.7) \\
		\hline
		\hline
		SCoPeS-SB: & $7.27 \cdot 10^1$ (4090) & $2.57 \cdot 10^1$ (3265) & $4.92 \cdot 10^1$ (2733) \\
		\hline
		\hline
	\end{tabular}
	\caption{SCoPeS-D results, Darcy approximation for samples S2\_100\_60, S3\_100\_58 and S3\_90\_53 of Category A. Permeability $k_{\mathrm{zz}}^{eff}$ in $mkDa$ (CPU time in s), $rtol_S = 10^{-9}$, $L=0.0009$ m, nproc=8. BC: pressure drop 1 Pa. The last line for comparison recalls permeability and CPU time when solving Stokes-Brinkman equations.}
	\label{table:Darcy-S2S3-A}
\end{table}

\noindent{\bf Solving Stokes-Brinkman equations for ternary images of Category B (Stokes connectivity).}
Let us now discuss the simulation results for samples from Category B (Stokes connectivity). Consider remaining images obtained from samples {\it S2} and {\it S3}. Again, the Stokes-Brinkman equations are solved for them, and the computed effective permeability and the CPU time are reported in Table \ref{table:S2S3-B}. 
One can observe that for these low porosity samples, SCoPeS-SB demonstrates very good performance. The computed effective permeability values are close to those computed with GeoDict. The convergence with respect to decreasing the tolerance is pronounced and stable. The increasing of the computational time when decreasing the tolerance is very moderate.

\begin{table}[h!]
  \centering
  \begin{tabular}{||c c c r | c c ||}
    \hline
    \multicolumn{6}{|| c ||}{Sample S2\_90\_55, Perm of porous voxels 236.3 mkDa} \\
    \hline
    \multicolumn{4}{|| c |}{GeoDict} & \multicolumn{2}{ c ||}{SCoPeS-SB} \\
    \hline
    Tol & Solver & $k_{\mathrm{zz}}^{eff}$, $mkDa$ & (CPU, s) & $rtol_S$ & $k_{\mathrm{zz}}^{eff}$, $mkDa$ (CPU, s) \\ [0.5ex]
    \hline
    \hline
    $ 10^{-1}$ & SimpleFFT: & $6.40\cdot10^{2}$ & 13996 & $ 10^{-5}$ & $5.26\cdot10^{2}$ (1853) \\
    \hline
    $ 1.1 \cdot 10^{-2}$ & SimpleFFT: & $5.95\cdot10^{2}$ & 227036 & $ 10^{-6}$ & $5.93\cdot10^{2}$ (2343) \\
    \hline
     & & & & $ 10^{-8}$ & $5.91\cdot10^{2}$ (3520) \\
    \hline
     & & & & $ 10^{-9}$ & $5.91\cdot10^{2}$ (4243) \\
    \hline
    \hline
    \hline
    \multicolumn{6}{|| c ||}{Sample S2\_80\_49, Perm of porous voxels 97.8 mkDa } \\
    \hline
    \multicolumn{4}{|| c |}{GeoDict} & \multicolumn{2}{ c ||}{SCoPeS-SB} \\
    \hline
    Tol & Solver: & $k_{\mathrm{zz}}^{eff}$, $mkDa$ & (CPU, s) & $rtol_S$ & $k_{\mathrm{zz}}^{eff}$, $mkDa$ (CPU, s) \\ [0.5ex]
    \hline
    \hline
    $10^{-1}$ & SimpleFFT: & $6.13\cdot10^{3}$ & 2103 & $ 10^{-5}$ & $5.95\cdot10^{3}$ (1456) \\
    \hline
    $ 10^{-2}$ & SimpleFFT: & $5.88\cdot10^{3}$ & 7347 & $ 10^{-6}$ & $5.81\cdot10^{3}$ (1816) \\
    \hline
    $ 10^{-3}$ & SimpleFFT: & $5.84\cdot10^{3}$ & 35744 & $10^{-7}$ & $5.83\cdot10^{3}$ (2319) \\
    \hline
    $10^{-1}$ & LIR: & $6.07\cdot10^{3}$ & 10509 & $10^{-8}$ & $5.83\cdot10^{3}$ (2892) \\
    \hline
    $10^{-2}$ & LIR: & $5.86\cdot10^{3}$ & 33216 & & \\
    \hline
    $10^{-3}$ & LIR: & $5.81\cdot10^{3}$ & 241828 & & \\
    \hline
    \hline
    \hline
    \multicolumn{6}{|| c ||}{Sample S3\_80\_48, Perm of porous voxels 84.4 mkDa } \\
    \hline
    \multicolumn{4}{|| c |}{GeoDict} & \multicolumn{2}{ c ||}{SCoPeS-SB} \\
    \hline
    Tol & Solver: & $k_{\mathrm{zz}}^{eff}$, $mkDa$ & (CPU, s) & $rtol_S$ & $k_{\mathrm{zz}}^{eff}$, $mkDa$ (CPU, s) \\ [0.5ex]
    \hline
    \hline
    $10^{-1}$ & SimpleFFT: & $1.89\cdot10^{4}$ & 568 & $10^{-3}$ & $1.63\cdot10^{4}$ (703) \\
    \hline
    $10^{-2}$ & SimpleFFT: & $1.56\cdot10^{4}$ & 5912 & $10^{-4}$ & $1.60\cdot10^{4}$ (955) \\
    \hline
    $10^{-3}$ & SimpleFFT: & $1.55\cdot10^{4}$ & 14099 & $10^{-5}$ & $1.55\cdot10^{4}$ (1190) \\
    \hline
     & & & & $10^{-7}$ & $1.55\cdot10^{4}$ (1945) \\
    \hline
    \hline
  \end{tabular}
  \caption{Ternary samples S2 and S3 of Category B. Permeability $k_{\mathrm{zz}}^{eff}$ in $mkDa$ computed with GeoDict (solvers SimpleFFT and LIR with periodic bc), and with SCoPeS-SB with pressure drop bc, $L=0.0009$ m, nproc=8.}
  \label{table:S2S3-B}
\end{table}


\noindent{\bf Solving Darcy approximation for ternary images of Category B (Stokes connectivity).}
For comparison, the Darcy approximation is also used for computing the effective permeability of samples from Category B. The results are summarized in Table \ref{table:Darcy-S2S3-B}. As expected, the Darcy approximation is not applicable for computing the effective permeability of the considered samples, namely samples for which Stokes connectivity exists. This illustrates the importance of the image classification stage in the presented workflow.


\begin{table}[h!]
  \centering
  \begin{tabular}{||c c c c||}
    \hline
    \hline
    $K_{Stokes}$, $mkDa$ & S2\_90\_55 & S2\_80\_49 & S3\_80\_48 \\ [0.5ex]
    \hline
    \hline
    $10^5$ & $6.13\cdot10^{1}$ (43.3) & $8.19\cdot10^{1}$ (45.0) & $1.39\cdot10^{2}$ (42.5) \\
    \hline
    $10^7$ & $1.56\cdot10^{3}$ (43.4) & $5.61\cdot10^{3}$ (43.1) & $1.05\cdot10^{4}$ (42.7) \\
    \hline
    $10^9$ & $1.45\cdot10^{5}$ (43.6) & $5.57\cdot10^{5}$ (42.9) & $1.05\cdot10^{6}$ (43.4) \\
    \hline
    $10^{11}$ & $1.45\cdot10^{7}$ (44.5) & $5.57\cdot10^{7}$ (43.6) & $1.05\cdot10^{8}$ (43.4) \\
    \hline
    \hline
    {SCoPeS-SB:} & $5.91\cdot10^{2}$ (4243) & $5.83\cdot10^{3}$ (3600) & $1.55\cdot10^{4}$ (1945) \\
    \hline
    \hline
  \end{tabular}
  \caption{SCoPeS-D results, Darcy approximation for samples S2\_90\_55, S2\_80\_49 and S3\_80\_48 of Category B. Permeability $k_{\mathrm{zz}}^{eff}$ in $mkDa$ (CPU time in s), $rtol_S=10^{-9}$, $L=0.0009$ m, nproc=8. BC: pressure drop 1 Pa. The last line for comparison recalls permeability and CPU time when solving Stokes-Brinkman equations.}
  \label{table:Darcy-S2S3-B}
\end{table}

\subsubsection{Simulations on multiclass image of size $300^3$.}
\label{subsubsec:multiclass}

The developed solver is able to work directly on multiclass images, such that they have individual permeability values in each unresolved voxel. On the one hand, this imposes higher memory requirements to the solver, on the other hand, this feature of the solver might be essential for certain classes of rocks. To illustrate this option, simulations are performed with image {\it S2}, as a multiclass image having 100 classes of porosity to account for the sub-micron scale effects. The results are summarized in Table \ref{table:S2}. Robust convergence of the simulations with respect to the tolerance in the case of multiclass image can be observed. The effective permeabilities computed on the three ternary images produced from the multiclass image {\it S2} are also provided in the table. It can be seen that in this particular case the effective permeabilities computed on the multiclass and on the three different ternary images differ significantly. This means that results obtained from simulations on ternary images might not be a good approximation of the results obtained on the true, multiclass grey image, at least if a simple averaging is used when producing ternary image out of a multiclass image.  

\begin{table}[h!]
  \centering
  \begin{tabular}{||c c c||}
    \hline
    \hline
    $rtol_S$ & $k_{\mathrm{zz}}^{eff}$, $mkDa$ & CPU time, s \\ [0.5ex]
    \hline
    \hline
    $ 10^{-5}$ & $8.63\cdot10^{2}$ & 1421 \\
    \hline
    $ 10^{-6}$ & $5.32\cdot10^{2}$ & 1883 \\
    \hline
    $ 10^{-7}$ & $5.29\cdot10^{2}$ & 2564 \\
    \hline
    $ 10^{-8}$ & $5.28\cdot10^{2}$ & 3444 \\
    \hline
    $ 10^{-9}$ & $5.28\cdot10^{2}$ & 4285 \\
    \hline
    \hline
    S2\_100\_60, $ 10^{-8}$ & $7.27\cdot10^{1}$ & 3499 \\
    \hline
    S2\_90\_55, $ 10^{-8}$ & $5.91\cdot10^{2}$ & 3520 \\
    \hline
    S2\_80\_49, $ 10^{-7}$ & $5.83\cdot10^{3}$ & 2319 \\
    \hline
  \end{tabular}
  \caption{SCoPeS-SB results, permeability $k_{\mathrm{zz}}^{eff}$ in $mkDa$ for multiclass sample S2, $L = 0.0009$ m, nproc = 8. BC: pressure drop 1 Pa. Permeabilities of the corresponding ternary samples are also reminded.}
  \label{table:S2}
\end{table}

\subsubsection{Simulations on ternary images of size $600^3$ and $1350^3$.}
\label{subsubsec:ternary_600_1350}

The simulations with larger images show the same behaviour as those with $300^3$ images. The simulation results for these images are presented in the Appendix. The classification of the $600^3$ with respect to flow connectivity is presented in \ref{table:600connectivity}
The simulation results for images  {\it T1,T2} of size $600^3$ are summarized in Table \ref{table:T1T2-A} for Category A and in Table \ref{table:T1T2-A} for Category B. Similarly to the simulation results for samples of size $300^3$, one can observe a robust convergence of the SCoPeS-SB solver with respect to the tolerance decrease. It also can be observed that the computational time is still relatively low, especially compared to GeoDict SimpleFFT solver. GeoDict-LIR solver is not explored here, we expect similar performance as for $300^3$ images. Comparing computational times reported, e.g., in Tables \ref{table:T1T2-A} and \ref{table:S1-A}, one can see that the computational time has increased about six times, what is very good result, having in mind the increase of the size of the samples. 
In general, when using AMG preconditioner, one can expect the computation time to increase proportionally to the number of unknowns. Thus, eight times increase of the CPU time was expected here. 
The fact that the increase is only six times can be explained by the fact that the more computationally intensive task exhibits better parallel scalability.
Simulations with Darcy approximation of Stokes-Brinkman model in the case of $600^3$ images show similar performance as in the case of $300^3$ images. Results of Darcy simulations are collected in Table \ref{table:Darcy-T1T2-A} for Category A and in Table \ref{table:Darcy-T1T2-B} for Category B.

Similarly, the results for {\it U1} sample of size $1350^3$ (Category A) are summarized in Table \ref{table:U1-A} for the Stokes-Brinkman problem and in Table \ref{table:Darcy-U1-A} for the Darcy problem.









\subsection{Summary.}
\label{subsec:summary}

A workflow containing image classification stage and customized efficient solvers for simulations on 3D CT images of samples from tight reservoirs with multiscale porosity is developed and intensively tested. The proposed stage for classifying images with respect to the existence of Stokes percolation patches allows for selecting the proper solver for each image. The developed solvers are validated using data from the literature and compared with a commercial software tool. Extensive testing on samples from a real tight reservoir are performed. The presented Stokes-Brinkman solver shows very robust convergence and high efficiency. This solver enables simulations not only on ternary, but also on multiclass images. The latter are directly extracted from grey images, and have individual permeability in each porous voxel. 
It was shown that in the case of missing pure fluid percolation, the Darcy approximation of the Stokes-Brinkman problem can be solved. This gives 10-15 times acceleration of the computations preserving the accuracy. 

\bibliographystyle{spmpsci}      
\bibliography{main.bib}   


\newpage

\setcounter{table}{0}
\setcounter{figure}{0}
\renewcommand{\thetable}{A.\arabic{table}}
\renewcommand{\thefigure}{A.\arabic{figure}}
\appendix

\noindent
{\Large Appendix} \\

This Appendix summarizes some further simulation results.

\section{Simulations on S1 ternary images of size $300^3$.}

All three samples produced from the S1 image belong to category A, no Stokes connectivity. Thus it is expected that the results in this case will be similar to those for S2\_100\_60, S2\_90\_55 and S3\_90\_53, which also belong to category A. Indeed, the results from Table \ref{table:S1-A} below are similar to those from Table \ref{table:SamplesS123} from the main text. 

\begin{table}[h!]
  \centering
  \begin{tabular}{||c c c r | c c ||}
    \hline
    \multicolumn{6}{|| c ||}{Sample S1\_100\_61, Perm of porous voxels 571.2 mkDa} \\
    \hline
    \multicolumn{4}{|| c |}{GeoDict} & \multicolumn{2}{ c ||}{SCoPeS-SB} \\
    \hline
    Tol & Solver & $k_{\mathrm{zz}}^{eff}$, $mkDa$ & CPU, s & $rtol_S$ & $k_{\mathrm{zz}}^{eff}$, $mkDa$ (CPU, s) \\ [0.5ex]
    \hline
    \hline
    $ 10^{-1}$ & SimpleFFT: & $8.74\cdot10^{1}$ & 4471 & $ 10^{-6}$ & $3.42\cdot10^{1}$ (1210) \\
    \hline
    $ 10^{-2}$ & SimpleFFT: & $6.16\cdot10^{1}$ & 61406 & $ 10^{-7}$ & $5.69\cdot10^{1}$ (1522) \\
    \hline
    $ 10^{-1}$ & LIR & $1.24\cdot10^{2}$ & 33968 & $ 10^{-8}$ & $5.68\cdot10^{1}$ (1898) \\
    \hline
    $ 9.95 \cdot 10^{-2}$ & LIR & $6.46\cdot10^{1}$ & 459561 & $ 10^{-9}$ & $5.68\cdot10^{1}$ (2390) \\
    \hline
    \hline
    \hline
    \multicolumn{6}{|| c ||}{Sample S1\_90\_56, Perm of porous voxels 273.8 mkDa } \\
    \hline
    \multicolumn{4}{|| c |}{GeoDict} & \multicolumn{2}{ c ||}{SCoPeS-SB} \\
    \hline
    Tol & Solver: & $k_{\mathrm{zz}}^{eff}$, $mkDa$ & CPU, s & $rtol_S$ & $k_{\mathrm{zz}}^{eff}$, $mkDa$ (CPU, s) \\ [0.5ex]
    \hline
    \hline
    $ 10^{-1}$ & SimpleFFT: & $1.17\cdot10^{2}$ & 5357 & $ 10^{-6}$ & $5.82\cdot10^{1}$ (1227) \\
    \hline
    $ 2 \cdot 10^{-2}$ & SimpleFFT: & $7.11\cdot10^{1}$ & 70690 & $ 10^{-7}$ & $6.45\cdot10^{1}$ (1584) \\
    \hline
     & & & & $ 10^{-8}$ & $6.34\cdot10^{1}$ (1942) \\
    \hline
     & & & & $ 10^{-9}$ & $6.33\cdot10^{1}$ (2431) \\
    \hline
    \hline
    \hline
    \multicolumn{6}{|| c ||}{Sample S1\_80\_50, Perm of porous voxels 113.3 mkDa } \\
    \hline
    \multicolumn{4}{|| c |}{GeoDict} & \multicolumn{2}{ c ||}{SCoPeS-SB} \\
    \hline
    Tol & Solver: & $k_{\mathrm{zz}}^{eff}$, $mkDa$ & CPU, s & $rtol_S$ & $k_{\mathrm{zz}}^{eff}$, $mkDa$ (CPU, s) \\ [0.5ex]
    \hline
    \hline
    $ 10^{-1}$ & SimpleFFT: & $8.07\cdot10^{1}$ & 16532 & $ 10^{-6}$ & $6.42\cdot10^{1}$ (1292) \\
    \hline
    $ 5 \cdot 10^{-3}$ & SimpleFFT: & $6.51\cdot10^{1}$ & 228113 & $ 10^{-7}$ & $6.03\cdot10^{1}$ (1628) \\
    \hline
     & & & & $ 10^{-8}$ & $7.01\cdot10^{1}$ (1991) \\
    \hline
     & & & & $ 10^{-9}$ & $6.00\cdot10^{1}$ (2393) \\
    \hline
    \hline
  \end{tabular}
  \caption{Ternary samples S1 of Category A. Permeability $k_{\mathrm{zz}}^{eff}$ in $mkDa$ computed with GeoDict (solvers SimpleFFT and LIR with periodic bc) and with SCoPeS-SB with pressure drop bc, $L=0.0009$ m, nproc=8.}
  \label{table:S1-A}
\end{table}


As mentioned above, all three samples produced from the S1 image belong to Category A. Thus, it is expected that the accuracy and the performance of the Darcy approximation in this case will be similar to those for S2\_100\_60, S2\_90\_55, and S3\_90\_53, which also belong to category A. Indeed, the results from Table \ref{table:Darcy-S1-A} below are similar to those from Table \ref{table:Darcy-S2S3-A} from the main text. 

\begin{table}[h!]
	\centering
	\begin{tabular}{||c c c c||}
		\hline
		\hline
		$K_{Stokes}$, $mkDa$ & S1\_100\_61 & S1\_90\_56 & S1\_80\_50 \\ [0.5ex]
		\hline
		\hline
		$10^5$ & $4.54 \cdot 10^1$ (43.2) & $4.48 \cdot 10^1$ (43.5) & $4.07 \cdot 10^1$ (42.6) \\
		\hline
		$ 10^7$ & $5.70 \cdot 10^1$ (43.6) & $6.35 \cdot 10^1$ (43.6) & $6.01 \cdot 10^1$ (43.3) \\
		\hline
		$10^9$ & $5.72 \cdot 10^1$ (44.1) & $6.39 \cdot 10^1$ (44.1) & $6.04 \cdot 10^1$ (43.6) \\
		\hline
		$10^{10}$ & $5.72 \cdot 10^1$ (44.5) & $6.39 \cdot 10^1$ (43.7) & $6.04 \cdot 10^1$ (43.1) \\
		\hline
		\hline
		SCoPeS-SB: & $5.68 \cdot 10^1$ (2390) & $6.34 \cdot 10^1$ (2431) & $6.00 \cdot 10^1$ (2393) \\
		\hline
		\hline
	\end{tabular}
	\caption{SCoPeS-D results, Darcy approximation for samples S1\_100\_61, S1\_90\_56 and S1\_80\_50 of Category A. Permeability $k_{\mathrm{zz}}^{eff}$ in $mkDa$ (CPU time in s), $rtol_S = 10^{-9}$, $L=0.0009$ m,  nproc=8. BC: pressure drop 1 Pa. The last line for comparison recalls permeability and CPU time when solving Stokes-Brinkman equations.}
	\label{table:Darcy-S1-A}
\end{table}

\newpage
\section{Darcy approximation, sensitivity study}

As mentioned, Darcy approximation of Stokes-Brinkman equations is used for images of Category A. Its accuracy and performance were examined before its further usage. 
Simulation results from sensitivity study on how the artificial permeability in the Darcy approximation influences the accuracy of the computations are presented in Table \ref{table:Darcy_S1} in the Appendix. The sample S1\_100\_61 is considered here.
This Table contains three subtables. In the first one, the fictitious permeability $K_{Stokes}$ in the Stokes (pure fluid) voxels is fixed to a moderate value $10^{-6}$ $m^2$, and the tolerance for the iterative method is varied. In the second subtable, a relatively rough tolerance is fixed for the iterative method ($10^{-4}$) and the fictitious permeability $K_{Stokes}$ is varied. Finally, in the third subtable, the tolerance $10^{-9}$ is fixed for the iterative method, and the fictitious permeability is varied. Comparing to the results from Table \ref{table:Darcy_S1}, one can see that Darcy approximation can be successfully used for fast and accurate computation of the effective permeability of samples in the case of no Stokes connectivity. The results from Table \ref{table:Darcy_S1} show that relatively low value for the fictitious permeability should be set in the Stokes voxels, and the iterative method should be converged with high accuracy. The latter however, in this particular case does not influence essentially the computational time.

\begin{table}[h!]
	\centering
	\begin{tabular}{||c c c c||}
		\hline
		\hline
		$rtol_S$ & $K_{Stokes}$, $mkDa$ & $k_{\mathrm{zz}}^{eff}$, $mkDa$ & CPU time, s \\ [0.5ex]
		\hline
		\hline
		$ 10^{-4}$ & $ 10^{6}$ & $5.51 \cdot 10^{1}$ & {39.2} \\
		\hline
		$ 10^{-5}$ & $ 10^{6}$ & $5.51 \cdot 10^{1}$ & 40.2 \\
		\hline
		$ 10^{-6}$ & $ 10^{6}$ & $5.51 \cdot 10^{1}$ & 40.6 \\
		\hline
		\hline
		\hline
		$ 10^{-4}$ & $ 10^{6}$ & {$5.51 \cdot 10^{1}$} & 39.2 \\
		\hline
		\hline
		$ 10^{-4}$ & $ 10^{8}$ & $5.75 \cdot 10^{1}$ & 38.8 \\
		\hline
		\hline
		$ 10^{-4}$ & $ 10^{10}$ & $2.31 \cdot 10^{2}$ & 39.3 \\
		\hline
		\hline
		\hline
		$ 10^{-9}$ & $ 10^{5}$ & $4.54 \cdot 10^{1}$ & 43.2 \\
		\hline
		$ 10^{-9}$ & $ 10^{6}$ & $5.51 \cdot 10^{1}$ & 43.4 \\
		\hline
		$ 10^{-9}$ & $ 10^{11}$ & $5.72 \cdot 10^{1}$ & 44.4 \\
		\hline
		\hline
	\end{tabular}
	\label{table:Darcy-S1-parameter_study}
	\caption{SCoPeS-D results, Darcy approximation for sample S1\_100\_61; $L = 0.0009$ m,  nproc=8. BC: pressure drop 1 Pa.}
	\label{table:Darcy_S1}
\end{table}


\section{Visual comparison of flow simulations with Stokes-Brinkman system and with its Darcy approximation}

\begin{figure}
     \centering
     \begin{subfigure}[b]{0.49\textwidth}
         \centering
         \includegraphics[width=\textwidth]{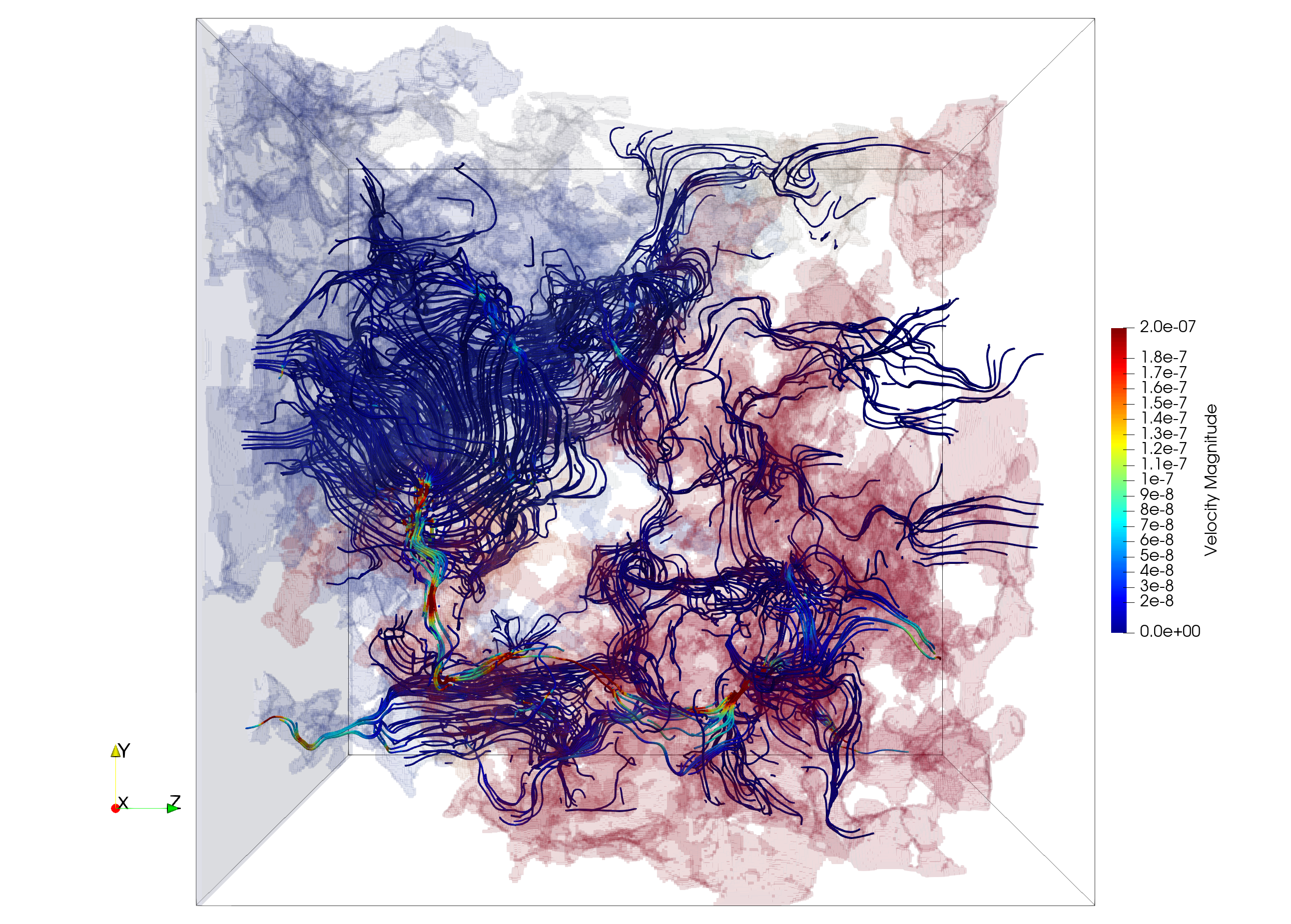}
         \caption{sample S2\_100\_60}
         \label{fig:s1}
     \end{subfigure}
     \begin{subfigure}[b]{0.49\textwidth}
         \centering
         \includegraphics[width=\textwidth]{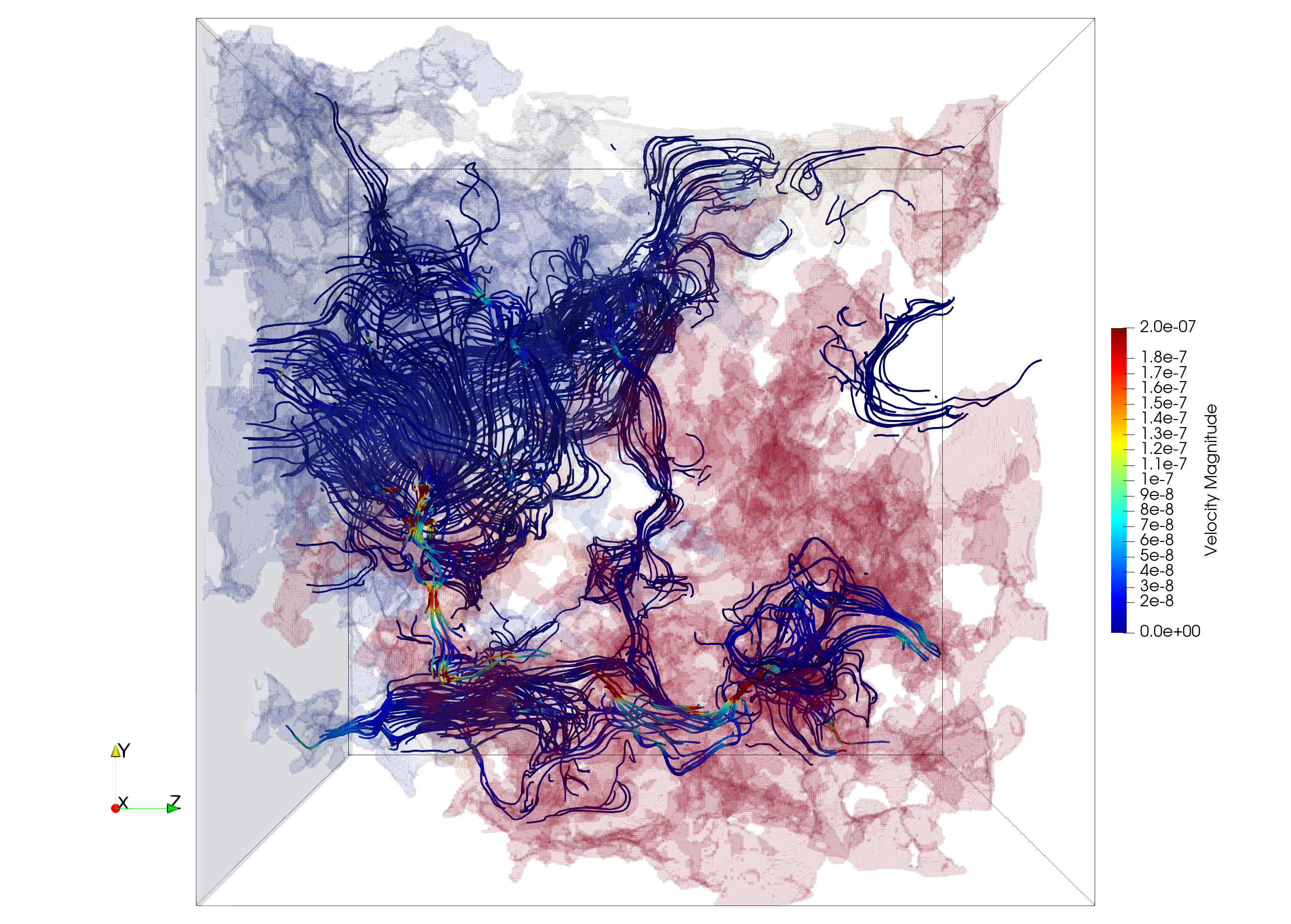}
         \caption{sample S2\_100\_60 darcy approximation}
         \label{fig:s2}
     \end{subfigure}
        \caption{sample S2\_100\_60}
        \label{fig:SBvsD}
\end{figure}

It was shown above that Darcy approximation can be efficiently used to compute integral quantities, such as permeability, for samples with no Stokes connectivity. Streamlines from computations with Stokes-Brinkman system and with its Darcy approximation are plotted on Fig. \ref{fig:SBvsD}. It can be seen that Darcy approximation not only accurately compute integral quantities, but also handles well the local behaviour of the flow. 

\section{Hardware/Software Specifications}

\begin{table}[h!]
  \centering
  \begin{tabular}{||c | c ||}
    \hline
    \multicolumn{2}{|| c ||}{Samples $300^3,600^3$} \\
    \hline
    Hardware model & Operating System \\
    \hline
    \begin{tabular}{@{}c@{}}2 x Intel Xeon E5-2687W v4@3.00 GHz, \\ RAM 528 Gb DDR4 2400 Hz \end{tabular}
     & Linux 4.15.0-154-generic Ubuntu 18.04.6 \\
    \hline
    \hline
    \multicolumn{2}{|| c ||}{Samples $1350^3$} \\
    \hline
    \begin{tabular}{@{}c@{}} Dell PowerEdge R640, \\ 2 x Intel Xeon Gold 6150@2.7 GHz, \\ RAM 1536 Gb DDR4 \end{tabular}
     & CentOS Linux release 7.9.2009 \\
    \hline
    \hline
  \end{tabular}
  \caption{Hardware/Software Specification}
  \label{table:HardSoft_spec}
\end{table}

\section{Simulations on ternary images of size $600^3$ and $1350^3$}

\begin{table}[h!]
	\centering
	\begin{tabular}{||c c c c c||} 
		\hline
		\hline
		Samples & $\phi$ (Res, Unres) & Stokes  & Stokes-Brinkman & File name \\ [0.5ex] 
	  &   &   Connectivity &   Connectivity & \\ [0.5ex]		\hline
		\hline
		{T1\_100\_59} & 0.141 (0.051, 0.089) & No & Yes & T1\_100\_59.raw \\
		\hline
		{T1\_90\_54} & 0.141 (0.062, 0.079) & No & Yes & T1\_90\_54.raw \\
		\hline
		{\bf T1\_80\_49} & 0.141 (0.073, 0.067) & {\bf Yes} & Yes & T1\_80\_49.raw \\
		\hline
		{T2\_100\_58} & 0.142 (0.051, 0.091) & No & Yes & T2\_100\_58.raw \\
		\hline
		{T2\_90\_54} & 0.142 (0.061, 0.081) & No & Yes & T2\_90\_54.raw \\
		\hline
		{\bf T2\_80\_49} & 0.142 (0.072, 0.070) & {\bf Yes} & Yes & T2\_80\_49.raw \\
		\hline
		\hline
		{U1\_100\_59} & 0.143 (0.052, 0.091) & No & Yes & U1\_100\_59.raw \\
		\hline
	\end{tabular}
	\caption{Porosity and Connectivity of the samples}
	\label{table:600connectivity}
\end{table}

\begin{table}[h!]
 \centering
 \begin{tabular}{||c c c r | c c ||}
  \hline
  \multicolumn{6}{|| c ||}{Sample T1\_100\_59, Perm of porous voxels 425.6 $mkDa$} \\
  \hline
  \multicolumn{4}{|| c |}{GeoDict} & \multicolumn{2}{ c ||}{SCoPeS-SB} \\
  \hline
  Tol & Solver & $k_{\mathrm{zz}}^{eff}$, $mkDa$ & (CPU, s) & $rtol_S$ & $k_{\mathrm{zz}}^{eff}$, $mkDa$ (CPU, s) \\ [0.5ex]
  \hline
  \hline
  $2.2 \cdot 10^{-1}$ & SimpleFFT & $6.49\cdot10^{1}$ & 170969 & $ 10^{-6}$ & $5.61\cdot10^{1}$ (6082) \\
  \hline
   & & & & $ 10^{-7}$ & $3.33\cdot10^{1}$ (8148) \\
  \hline
   & & & & $ 10^{-8}$ & $3.37\cdot10^{1}$ (10538) \\
  \hline
   & & & & $ 10^{-9}$ & $3.38\cdot10^{1}$ (13027) \\
  \hline
  \hline
  \multicolumn{6}{|| c ||}{Sample T1\_90\_54, Perm of porous voxels 204.0 $mkDa$} \\
  \hline
  \multicolumn{4}{|| c |}{GeoDict} & \multicolumn{2}{ c ||}{SCoPeS-SB} \\
  \hline
  Tol & Solver: & $k_{\mathrm{zz}}^{eff}$, $mkDa$ & (CPU, s) & $rtol_S$ & $k_{\mathrm{zz}}^{eff}$, $mkDa$ (CPU, s) \\ [0.5ex]
  \hline
  \hline
  $ 7.2 \cdot 10^{-1}$ & SimpleFFT: & $1.41\cdot10^{2}$ & 52369 & $10^{-6}$ & $3.38\cdot10^{1}$ (6059) \\
  \hline
   & & & & $ 10^{-7}$ & $5.51\cdot10^{1}$ (8035) \\
  \hline
   & & & & $ 10^{-8}$ & $5.48\cdot10^{1}$ (10240) \\
  \hline
  \hline
  \hline
  \multicolumn{6}{|| c ||}{Sample T2\_100\_58, Perm of porous voxels 367.4 $mkDa$} \\
  \hline
  \multicolumn{4}{|| c |}{GeoDict} & \multicolumn{2}{ c ||}{SCoPeS-SB} \\
  \hline
  Tol & Solver: & $k_{\mathrm{zz}}^{eff}$ & (CPU, s) & $rtol_S$ & $k_{\mathrm{zz}}^{eff}$, $mkDa$ (CPU, s) \\ [0.5ex]
  \hline
  \hline
  $ 2 \cdot 10^{-1}$ & SimpleFFT: & $4.25\cdot10^{1}$ & 1.2509e+6 & $10^{-6}$ & $3.61\cdot10^{1}$ (6191) \\
  \hline
   & & & & $ 10^{-7}$ & $3.35\cdot10^{1}$ (8320) \\
  \hline
   & & & & $ 10^{-8}$ & $3.36\cdot10^{1}$ (10425) \\
  \hline
  \hline
  \hline
  \multicolumn{6}{|| c ||}{Sample T2\_90\_54, Perm of porous voxels 204.0 $mkDa$} \\
  \hline
  \multicolumn{4}{|| c |}{GeoDict} & \multicolumn{2}{ c ||}{SCoPeS-SB} \\
  \hline
  Tol & Solver: & $k_{\mathrm{zz}}^{eff}$, $mkDa$ & (CPU, s) & $rtol_S$ & $k_{\mathrm{zz}}^{eff}$, $mkDa$ (CPU, s) \\ [0.5ex]
  \hline
  \hline
  $ 2.4 \cdot 10^{-1}$ & SimpleFFT: & $9.32\cdot10^{1}$ & 815269 & $10^{-6}$ & $5.94\cdot10^{1}$ (6333) \\
  \hline
   & & & & $ 10^{-7}$ & $6.38\cdot10^{1}$ (8371) \\
  \hline
   & & & & $ 10^{-8}$ & $6.51\cdot10^{1}$ (10456) \\
  \hline
   & & & & $ 10^{-9}$ & $6.51\cdot10^{1}$ (13191) \\
  \hline
  \hline
 \end{tabular}
 \caption{Ternary samples T1 and T2 with no Stokes connectivity. Permeability $k_{\mathrm{zz}}^{eff}$ in $mkDa$ computed with GeoDict (solvers SimpleFFT and LIR with periodic bc) and with SCoPeS-SB with pressure drop bc, $L=0.0018$ m, nproc=8.}
 \label{table:T1T2-A}
\end{table}

\begin{table}[h!]
  \centering
  \begin{tabular}{||c c c c c||}
    \hline
    \hline
    $K_{Stokes}$, & T1\_100\_59 & T1\_90\_54 & T2\_100\_58 & T2\_90\_54 \\ [0.5ex]
    $mkDa$ &  &  &  &  \\
    \hline
    \hline
    $10^5$ & $2.58\cdot10^{1}$ (260) & $2.18\cdot10^{1}$ (319) & $2.17\cdot10^{1}$ (321) & $2.36\cdot10^{1}$ (322) \\
    \hline
    $ 10^7$ & $3.39\cdot10^{1}$ (255) & $5.84\cdot10^{1}$ (318) & $3.36\cdot10^{1}$ (330) & $7.23\cdot10^{1}$ (327) \\
    \hline
    $10^9$ & $3.40\cdot10^{1}$ (264) & $6.21\cdot10^{1}$ (325) & $3.39\cdot10^{1}$ (334) & $7.77\cdot10^{1}$ (332) \\
    \hline
    $10^{11}$ & $3.41\cdot10^{1}$ (269) & $6.22\cdot10^{1}$ (320) & $3.39\cdot10^{1}$ (334) & $7.77\cdot10^{1}$ (331) \\
    \hline
    \hline
    SCoPeS-SB: & $3.30\cdot10^{1}$ (9472) & $4.33\cdot10^{1}$ (10425) & $3.36\cdot10^{1}$ (10425) & $6.51\cdot10^{1}$ (10456) \\
    \hline
    \hline
  \end{tabular}
  \caption{SCoPeS-D results, Darcy approximation for samples T1\_100\_59, T1\_90\_54, T2\_100\_58 and T2\_90\_54. Permeability $k_{\mathrm{zz}}^{eff}$ in $mkDa$ (CPU time in s), $rtol_S = 10^{-9}$, $L=0.0018$ m, nproc=8. BC: pressure drop 1 Pa. The last line for comparison recalls permeability and CPU time when solving Stokes-Brinkman equations.}
  \label{table:Darcy-T1T2-A}
\end{table}

\newpage

\begin{table}[h!]
  \centering
  \begin{tabular}{||c c c r | c c ||}
    \hline
    \multicolumn{6}{|| c ||}{Sample T1\_80\_49, Perm of porous voxels 97.8 $mkDa$} \\
    \hline
    \multicolumn{4}{|| c |}{GeoDict} & \multicolumn{2}{ c ||}{SCoPeS-SB} \\
    \hline
    Tol & Solver & $k_{\mathrm{zz}}^{eff}$, $mkDa$ & (CPU, s) & $rtol_S$ & $k_{\mathrm{zz}}^{eff}$, $mkDa$ (CPU, s) \\ [0.5ex]
    \hline
    \hline
    $ 10^{-1}$ & SimpleFFT & $6.45\cdot10^{2}$ & 155062 & $ 10^{-5}$ & $5.68\cdot10^{2}$ (4849) \\
    \hline
    $1.3 \cdot 10^{-2}$ & SimpleFFT & $6.01\cdot10^{2}$ & 834989 & $ 10^{-6}$ & $6.06\cdot10^{2}$ (6208) \\
    \hline
     & & & & $ 10^{-7}$ & $5.88\cdot10^{2}$ (8215) \\
    \hline
     & & & & $ 10^{-8}$ & $5.88\cdot10^{2}$ (10460) \\
    \hline
    \hline
    \hline
    \multicolumn{6}{|| c ||}{Sample T2\_80\_49, Perm of porous voxels 97.8 $mkDa$} \\
    \hline
    \multicolumn{4}{|| c |}{GeoDict} & \multicolumn{2}{ c ||}{SCoPeS-SB} \\
    \hline
    Tol & Solver: & $k_{\mathrm{zz}}^{eff}$, $mkDa$ & (CPU, s) & $rtol_S$ & $k_{\mathrm{zz}}^{eff}$, $mkDa$ (CPU, s) \\ [0.5ex]
    \hline
    \hline
    $ 10^{-1}$ & SimpleFFT: & $6.59\cdot10^{3}$ & 6087 & $10^{-4}$ & $1.94\cdot10^{3}$ (4238) \\
    \hline
    $2.6 \cdot 10^{-2}$ & SimpleFFT & $1.59\cdot10^{3}$ & 682309 & $ 10^{-5}$ & $1.49\cdot10^{3}$ (5352) \\
    \hline
     & & & & $ 10^{-6}$ & $1.57\cdot10^{3}$ (6899) \\
    \hline
     & & & & $ 10^{-7}$ & $1.57\cdot10^{3}$ (8801) \\
    \hline
    \hline
  \end{tabular}
  \caption{Ternary samples T1 and T2 of Category B. Permeability $k_{\mathrm{zz}}^{eff}$ in $mkDa$ computed with GeoDict (solvers SimpleFFT and LIR with periodic bc) and with SCoPeS-SB with pressure drop bc, $L=0.0018$ m, nproc=8.}
  \label{table:T1T2-B}
\end{table}

\begin{table}[h!]
  \centering
  \begin{tabular}{||c c c||}
    \hline
    \hline
    $K_{Stokes}$, $mkDa$ & T1\_80\_49 & T2\_80\_49 \\ [0.5ex]
    \hline
    \hline
    $10^5$ & $2.64\cdot10^{1}$ (314) & $3.45\cdot10^{1}$ (320) \\
    \hline
    $10^7$ & $1.04\cdot10^{3}$ (320) & $1.65\cdot10^{3}$ (321) \\
    \hline
    $10^9$ & $1.02\cdot10^{5}$ (325) & $1.60\cdot10^{5}$ (322) \\
    \hline
    $10^{11}$ & $1.02\cdot10^{7}$ (323) & $1.60\cdot10^{7}$ (327) \\
    \hline
    \hline
    SCoPeS-SB: & $5.89\cdot10^{2}$ (8215) & $1.57\cdot10^{3}$ (6899) \\
    \hline
    \hline
  \end{tabular}
  \caption{SCoPeS-D results, Darcy approximation for samples S2\_90\_55, S2\_80\_49 and S3\_80\_48. Permeability $k_{\mathrm{zz}}^{eff}$ in $mkDa$ (CPU time in s), $rtol_S = 10^{-9}$, $L=0.0018$ m, nproc=8. BC: pressure drop 1 Pa. The last line for comparison recalls permeability and CPU time when solving Stokes-Brinkman equations.}
  \label{table:Darcy-T1T2-B}
\end{table}


\begin{table}[h!]
  \centering
  \begin{tabular}{||c c c r | c c ||}
    \hline
    \multicolumn{6}{|| c ||}{Sample U1\_100\_59, Perm of porous voxels 425.60 $mkDa$} \\
    \hline
    \multicolumn{4}{|| c |}{GeoDict} & \multicolumn{2}{ c ||}{SCoPeS-SB} \\
    \hline
    Tol & Solver & $k_{\mathrm{zz}}^{eff}$, $mkDa$ & (CPU, s) & $rtol_S$ & $k_{\mathrm{zz}}^{eff}$, $mkDa$ (CPU, s) \\ [0.5ex]
    \hline
    \hline
    $3.92 \cdot 10^{-1}$ & SimpleFFT & $9.17\cdot10^{1}$ & 1.11362e+06 & $ 10^{-6}$ & $3.34\cdot10^{1}$ (99841) \\
    \hline
    $10^{-1}$ (nproc=16) & SimpleFFT & $6.06\cdot10^{1}$ & 865126 & $ 10^{-7}$ & $3.29\cdot10^{1}$ (129588) \\
    \hline
    $10^{-1}$ & LIR & divergence & & $ 10^{-8}$ & $3.28\cdot10^{1}$ (165821) \\
    \hline
     & & & & $ 10^{-9}$ & $3.28\cdot10^{1}$ (205090) \\
    \hline
    \hline
    \hline
  \end{tabular}
  \caption{Ternary sample U1 with no Stokes connectivity. Permeability $k_{\mathrm{zz}}^{eff}$ in $mkDa$ computed with GeoDict (solvers SimpleFFT with periodic bc) and with SCoPeS-SB with pressure drop bc, $L=0.0027$ m, nproc=8.}
  \label{table:U1-A}
\end{table}

\begin{table}[h!]
  \centering
  \begin{tabular}{||c c||}
    \hline
    \hline
    $K_{Stokes}^{-1}$, $mkDa$ & U1\_100\_59 \\ [0.5ex]
    \hline
    \hline
    $10^5$ & $2.35\cdot10^{1}$ (7255) \\
    \hline
    $10^7$ & $3.34\cdot10^{1}$ (7373) \\
    \hline
    $10^9$ & $3.37\cdot10^{1}$ (7537) \\
    \hline
    $10^{11}$ & $3.37\cdot10^{1}$ (7437) \\
    \hline
    \hline
    SCoPeS-SB & $3.28\cdot10^{1}$ (205090) \\
    \hline
    \hline
  \end{tabular}
  \caption{SCoPeS-D results, Darcy approximation for sample U1. Permeability $k_{\mathrm{zz}}^{eff}$ in $mkDa$ (CPU time in s), $rtol_S = 10^{-9}$, $L=0.0009$ m, nproc=8. BC: pressure drop 1 Pa. The last line for comparison recalls permeability and CPU time when solving Stokes-Brinkman equations.}
  \label{table:Darcy-U1-A}
\end{table}

\begin{table}[h!]
  \centering
  \begin{tabular}{||c c c c c||}
    \hline
    \hline
    sample & \# non-solid voxels & SCoPeS-SB (Gb) & SCoPeS-D (Gb) & GeoDict SimpleFFT (Gb) \\ [0.5ex]
    \hline
    \hline
    S1\_100\_61 & 3.79 mln. & 19.7 & 13.8 & 2.2 \\
    \hline
    T1\_100\_59 & 30.5 mln. & 104.8 & 87.5 & 16.5 \\
    \hline
    U1\_100\_59 & 350.8 mln. & 1025.2 & 896.3 & 183.6 \\
    \hline
    \hline
  \end{tabular}
  \caption{Memory (RAM) usage of SCoPeS comparing with GeoDict-SimpleFFT.}
  \label{table:memory}
\end{table}

\end{document}